\documentclass[final,reqno]{siamltex}

\usepackage{amsmath}
\usepackage{amssymb}
\usepackage{graphicx}
\newtheorem{remark}{Remark}[section]
\newtheorem{thm}{Theorem}[section]

\newtheorem{lem}{Lemma}[section]

\def\div{\mbox{div}}
\def\grad{\mbox{grad}}

\def\dfrac{\displaystyle\frac}

\begin{document}

\title{Boundary Shape Control of Navier-Stokes Equations and
Geometrical Design Method for Blade's Surface in the
Impeller\thanks{Subsidized by NSFC Project
50136030,50306019,40375010,10471110,10471109.}}

\author{Kaitai Li\thanks{School of Sciences, Xi'an Jiaotong University,Xi'an, 710049,P.R.China, (ktli@xjtu.edu.cn).}
 \and Jian Su\thanks{School of Sciences, Xi'an Jiaotong
University,Xi'an, 710049,P.R.China.} \and Liquan Mei\thanks{School
of Sciences, Xi'an Jiaotong University,Xi'an, 710049,P.R.China.}}

\maketitle

\setcounter{page}{1}

%\begin{PII}
%S00000000000
%\end{PII}

%%%%%%%%%
\begin{abstract}
In this paper A Geometrical Design Method for Blade's surface $\Im$
in the impeller is provided here $\Im$ is a solution to  a coupling
system consisting of the well-known Navier-Stokes equations and a
four order elliptic boundary value problem . The coupling system is
used to describe the relations between solutions of Navier-Stokes
equations and the geometry of the domain occupied by fluids, and
also provides new theory and methods for optimal geometric design
 of the boundary of domain mentioned above. This coupling system is the
Eular-Lagrange equations of the optimal control problem which is
describing a new principle of the geometric design for the blade's
surface of an impeller. The control variable is the surface of the
blade and the state equations are Navier-Stokes equations with mixed
boundary conditions in the channel between two blades. The objective
functional depending on the geometry shape of blade's surface
describes the  dissipation energy of the flow and the power of  the
impeller.  First we prove the  existence of a solution of the
optimal control problem. Then we use a special coordinate system of
the Navier-Stokes equations to derive  the objective functional
  which depends on the surface $\Theta$ explicitly. We also show the
 weakly continuity of the solution of the Navier-Stokes
 equations with respect to the geometry shape of the blade's surface.
 \end{abstract}



\pagestyle{myheadings} \thispagestyle{plain} \markboth{Kaitai Li,
Jian Su, and Liquan Mei}{Boundary Shape Control of Navier-Stokes
Equations}

%%%%%%%%%%%%%%%%%%%%%%%%%%%%%%%%%%%%%%%%%%%%%%%%%%%%%%%%%%%%%%%%%%%%%
%
\section{Introduction}\label{sec-1}

Blade's shape design for the impeller is driven by the need of
improving performance and reliability. So far we have not found a
geometric design entirely from mathematical point of view. As it is
well known that the blade's surface is a part of the boundary of the
flow domain in the impeller. We can use the technique of the
boundary geometric control problem for the Navier-Stokes equations
to design the blade's shape. This idea is motivated by the classical
minimal surface which is to find a surface spanning on a closed
Jordanian curvilinear $C$ such that
\[
J(\Im)=Aug\inf\limits_{S\,\in\,\mathcal{F}}J(S)
\]
where $J(S)=\int\int\limits_S dS$ is the area of $S$.

In this paper we try to propose a principle for a fully mathematical
design of the surface of the blade in an impeller. This principle
models a general minimal surface by minimizing a functional proposed
by us. A key point in this modelling process is theoretical
rationality and the realization of our design procedure. Using a
tensor analysis technique we realize this procedure and obtain the
Eular-Lagrange equations for blade's surface which is a system
coupling an elliptic boundary value problem, the Navier-Stokes
equations and linearized Navier-Stokes equations, and prove the
existence of solution of the system coupling problem.

This paper is organized as follows. In  section \ref{sec-2} we give
the main results of this paper. In section \ref{sec-3}, we derive
the rotating Navier-Stokes equations in the channel in the impeller
with mixed boundary  condition under a new coordinate system. We
give the minimizing functional problem and derive the Euler-Lagrange
equations in section \ref{sec-4}. In section \ref{sec-5} another
model to design the blade of the impeller is given. In the last
section we prove the existence of the solution to the optimal
control  problem, including the existence of the solution  of
Navier-Stokes Equations with mixed boundary conditions and the
weakly continuous dependence relationships of the solution of
Navier-Stokes  equation with respect to the geometry shape of the
blade's  surface and so on.
%%%%%%%%%%%%%%%%%%%%%%%%%%%%%%%%%%%%%%%%%
\section{Main Results}\label{sec-2}

Suppose $(x^1,x^2)\in\,D\subset E^2 $($2D-Euclidian Space$). A
smooth mapping $\Theta(x^1,x^2)$ is the image of a surface. On the
other hand, suppose that $(r,\theta,z)$ is a polar cylindrical
coordinate system rotating with impeller's angular velocity $\omega$
.

%\begin{center}
%\includegraphics[bb=0 0 3.77cm 6.75cm]{point.bmp}
%\end{center}

$(\vec{e}_r,\vec{e}_\theta,\vec{k})$ are the
corresponding base vectors. z-axis is the rotating axis of the
impeller. N is the number of blade and
$\varepsilon=\pi/N$. The angle between two successive blades is $\frac{2\pi%
}{N}$. The flow passage of the impeller is bounded by $\partial
\Omega_\varepsilon=\Gamma_{in}\cup\Gamma_{out}\cup\Gamma_t\cup\Gamma_b\cup
S_+\cup S_{-}$.  The middle surface $S$ of the blade is defined as
the image $\vec{\Re}$ of the closure of a domain $D\subset R^2$
where $\vec{\Re}: D\rightarrow \Re^3$ is a smooth injective mapping
which can be expressed by that for any point $\vec{\Re}(D)\in S$
\begin{equation}\label{lab2.1}
\vec{\Re}(x)
=x^2\vec{e}_r+x^2\Theta(x^1,x^2)\vec{e}_\theta+x^1\vec{k}, \forall
x=(x^1,x^2)\in \overset{-}{D}, \end{equation}
where $\Theta \in
C^2(D,R)$ is a smooth function. $x=(x^1,x^2)$ is called a Gaussian
coordinate system on $S$. It is easy to prove that there exists a
family $S_\xi$ of surfaces with a single parameter to cover the domain $%
\Omega_\varepsilon$ defined by the mapping $D\rightarrow S_\xi=\{\vec{R}%
(x^1,x^2;\xi):\,\forall (x^1,x^2)\in D\}$:
\begin{equation}\label{lab2.2}
\vec{R}(x^1,x^2;\xi)=x^2\vec{e}_r+x^2(\varepsilon\xi+ \Theta(x^1,x^2))\vec{e}%
_\theta+x^1\vec{k},
\end{equation}
It is clear that the
metric tensor $a_{\alpha\beta}$ of $S_\xi$ is homogenous and
nonsingular independent of $\xi$, and is given as follows:
\begin{equation}\label{lab2.3}
a_{\alpha\beta}=\frac{\partial \vec{R}}{\partial
x^\alpha}\frac{\partial \vec{R}}{\partial x^\beta}=
\delta_{\alpha\beta}+r^2\Theta_\alpha\Theta_\beta,\quad a=\mbox{det}%
(a_{\alpha\beta})=1  +r^2(\Theta_1^2+\Theta_2^2)>0,
\end{equation}
From this we establish a curvilinear coordinate system $(x^1,x^2,\xi)$ in $%
\Re^3$ ,
\begin{equation}\label{lab2.4}
(r,\theta,z)\rightarrow (x^1,x^2,\xi):  x^1=z,\quad x^2=r,\quad
\xi=\varepsilon^{-1}(\theta-\Theta(x^1,x^2)),
\end{equation} that
maps the flow passage domain
\[
\Omega_\varepsilon=\{\vec{R}(x^1,x^2,\xi)=x^2\vec{e}_r+
x^2(\varepsilon\xi+\Theta(x^1,x^2))\vec{e}_\theta+x^1\vec{k},
\forall (x^1,x^2,\xi)\in \Omega \}
\]
into a fixed domain in $E^3$(3D Euclidian Space):
\[
\Omega=\{(x^1,x^2)\in D, -1\leq\xi\leq 1\}\quad \mbox{in}\quad \Re^3
\]
which is independent of Surface $S$ of the blade, and Jacobian
\[
J(\frac{\partial (r,\theta,z)}{\partial(x^1,x^2,\xi)})=\varepsilon,
\]
therefore the transformation is nonsingular.

Assume that $(x^{1^{\prime }},x^{2^{\prime }},x^{3^{\prime }})=(r,\theta,z)$%
,as well known that corresponding metric tensor of $\Re^3$ is
$(g_{1^{\prime }1^{\prime }}=1,g_{2^{\prime }2^{\prime
}}=r^2,g_{3^{\prime }3^{\prime }}=1,g_{i^{\prime }j^{\prime }}=0
\forall i^{\prime }\neq j^{\prime })$. According to rule of tensor
transformation under coordinate transformation we have following
calculation formulae
\[
g_{ij}=g_{i^{\prime }j^{\prime }}\frac{\partial x^{i^{\prime
}}}{\partial x^{i}} \frac{\partial x^{j^{\prime }}}{\partial x^{j}}.
\]
Substituting \eqref{lab2.3} into above formulae the metric tensor of
$E^3$ in new curvilinear coordinate system can be obtain

\begin{equation}\label{lab2.5}
g_{\alpha\beta}=a_{\alpha\beta},\quad
g_{3\beta}=g_{\beta3}=\varepsilon r^2\Theta_\beta,\quad
g_{33}=\varepsilon^2 r^2,\quad g=\mbox{det}(g_{ij})=\varepsilon^2
r^2.
\end{equation}
Through this paper we denote
$\Theta_\alpha=\frac{\partial\Theta}{\partial x^\alpha}$.  Its
contravariant components are given by
\begin{equation}\label{lab2.6}
g^{\alpha\beta}=\delta^{\alpha\beta},\,
g^{3\beta}=g^{\beta3}=-\varepsilon^{-1}\Theta_\beta,\quad
g^{33}=\varepsilon^{-2}r^{-2}(1+r^2|\nabla\Theta|^2),
\end{equation}
where $|\nabla\Theta|^2=\Theta_1^2+\Theta^2_2\,\,\mbox{and}\,\,
\Theta_\alpha=\displaystyle\frac{\partial\Theta}{\partial
x^\alpha}$.

\begin{enumerate}
  \item[$\clubsuit$] First Model:
  \begin{thm}\label{thm1}
Suppose the $\Theta$ is a blade's surface defined by \eqref{lab2.1}.
Then $\Theta$ is proposed as a solution of following elliptic
boundary value problem:
\begin{equation}\label{lab2.7}
\left\{%
\begin{array}{ll}
\frac{\partial^2}{\partial x^\lambda\partial x^\sigma}(K^{\alpha\beta\lambda%
\sigma}(w) \frac{\partial^2\Theta}{\partial x^\beta\partial
x^\alpha}) +
\frac{\partial^2}{\partial x^\lambda \partial x^\sigma}(r\widehat{\Phi}%
^{\lambda\sigma}(w,\Theta)) &  \\
\quad\quad-\frac{\partial}{\partial x^\lambda}(r\hat{\Phi}%
^\lambda(w,\Theta)) +\hat{\Phi}^0(w,\widehat{w})r=0,\quad \forall
(x^1,x^2)\in\,D\subset \,\Re^2 &  \\
\Theta=\Theta_0,\quad \frac{\partial \Theta}{\partial n}=\Theta_*,
\text{on} \partial D
\end{array}%
\right. \end{equation}
 combing Navier-Stokes equations and
linearized Navier-Stokes equations, where $(w,p)$ and $(\widehat{w},
\widehat{p})$ are solutions of compressible or incompressible
rotating Navier-Stokes equations \eqref{lab3.1} or \eqref{lab3.7}
and linearized Navier-Stokes equations (3.34) respectively and
\begin{equation}\label{lab2.8}
K^{\alpha\beta\lambda\sigma}(w,\Theta)=2\mu
r^3W^{\alpha\sigma}\delta^{\beta\lambda},\quad
W^{\alpha\beta}=\int^1_{-1}w^\alpha w^\beta d\xi,
\end{equation}
$\widetilde{\Phi}^0,\widetilde{\Phi}^\lambda, \widetilde{\Phi}%
^{\lambda\sigma}$ are defined by (4.11) respectively.
\end{thm}

Variational formulation associated with \eqref{lab2.7} is given by
\begin{equation}\label{lab2.9}
\left\{%
\begin{array}{ll}
\text{Find } \Theta\,\in V_\Gamma(D)=\{q|q\in\,H^2(D),
q|_{\Gamma_0}=\Theta_0,\,\frac{\partial q}{\partial
n}|_{\Gamma_0}=0\},& \\
\text{such that, } \quad\forall\eta\,\in\,H^2_0(D) &\\
\int\int_D[ (K^{\lambda\sigma\nu\mu}(w) \Theta_{\nu\mu}+r\hat{\Phi}%
^{\lambda\sigma}(w,\Theta)] \eta_{\lambda\sigma}dx &  \\
\quad+\int\int_D[r\hat{\Phi}^\lambda(w,\Theta)\eta_\lambda+ r\hat{\Phi}^0(w,%
\widehat{w},\Theta)\eta]dx=0 .
\end{array}%
\right. \end{equation}
  \item[$\clubsuit$] Second Model:
  \begin{thm}\label{thm2}
  Suppose the $\Theta$ is a blade's surface defined by
\eqref{lab2.1}. Then $\Theta$ is proposed as a solution of following
elliptic boundary value problem :
\begin{equation}\label{lab2.10}
\left\{%
\begin{array}{ll}
-\left(K_0(w)\widetilde{\Delta}\Theta+K^{\lambda\nu}(w,\Theta)\Theta_{\nu%
\lambda}\right)+F^{\nu\mu}(w)\Theta_\nu\Theta_\mu&\\
\quad+F^\lambda(w)\Theta_\lambda+F_0(w,\Theta)=0, &  \\
\Theta|_\gamma=\Theta_0,
\end{array}%
\right.
\end{equation}
where $K_0(w), K^{\lambda\nu}(w), F^{\nu\mu}(w), F^\lambda(w),
F_0(w,\Theta) $ are defined by (5.15).
  \end{thm}
  The variational formulation associated with (5.11) is given by
\begin{equation}\label{lab2.11}
\left\{%
\begin{array}{ll}
\text{Find } \Theta \in H_\gamma^1(D)=\{v|v\in\,H^1(D),v=\Theta^*\text{ on }%
\gamma=\partial D\}&\\
\text{such that, } \forall\eta\in\,H^1_0(D)&  \\
\int\int_D\left\{[\Psi_0(w,p,\Theta)\eta+\Psi^\lambda(w,p,\Theta)\eta_\lambda
-\mu
r^2 W^\sigma\frac{\partial (\Theta_\lambda\eta_\lambda)}{\partial x^\sigma}%
]r\omega\varepsilon r\right\}dx=0,
\end{array}%
\right.
\end{equation}
where
\begin{equation}\label{lab2.12}
\Psi^\lambda(w,p,\Theta)=\Psi_0^\lambda(w,p)
+\Psi_\nu^\lambda(w,p)\Theta_\nu+
\Psi_{\nu\mu}^\lambda(w,p)\Theta_\nu\Theta_\mu,
\end{equation}
and $\Psi_0(w,p,\Theta),\Psi_0^\lambda(w,p) ,\Psi_\nu^\lambda(w,p),+
\Psi_{\nu\mu}^\lambda(w,p)$ are defined by (5.10)(5.11).
\end{enumerate}
%%%%%%%%%%%%%%%%%%%
\section{Rotating Navier-Stokes Equations With Mixed Boundary
Conditions}\label{sec-3}

At first, we consider the three-dimensional rotating Navier-Stokes
equations in a frame rotating around the axis of a  rotating
impeller with an angular velocity $\omega$:
\begin{equation}\label{lab3.1}
\left\{%
\begin{array}{ll}
\displaystyle\frac{\partial\rho}{\partial t}+\mbox{div}(\rho w)=0, &  \\
\rho a=\nabla\sigma+f, &  \\
\rho c_v(\displaystyle\frac{\partial T}{\partial t}+w^j\nabla_j T)-\mbox{div}%
(\kappa\mbox{grad} T)+p\mbox{div} w-\Phi=h, &  \\
p=p(\rho,T),
\end{array}%
\right.
\end{equation}

%\begin{center}
%\includegraphics[bb=0 0 6cm 6cm]{yelun.bmp} \includegraphics[bb=0 0 6cm 6cm]{yepian.bmp}
%\end{center}

where $\rho$ is the density of the fluid, $w$ the velocity of the
fluid, $h$ the heat source, $T$ the temperature, $k$ the coefficient
of heat conductivity, $C_v$ specific heat at constant volume, and
$\mu$ viscosity. Furthermore, the strain rate tensor, stress
tensor,dissipative function and viscous tensor are given by
respectively:
\begin{equation}\label{lab3.2}
\left\{
\begin{array}{ll}
e_{ij}(w) & =\displaystyle\frac12(\nabla_i w_j+\nabla_j w_i);\quad
i,j=1,2,3,
\\
e^{ij}(w) & =g^{ik}g^{jm}e_{km}(w)=\displaystyle\frac 12(\nabla^i
w^j+\nabla^j w^i), \\
\sigma^{ij}(w,p) & =A^{ijkm} e_{km}(w),\quad
\Phi=A^{ijkm}e_{ij}(w)e_{ij}(w),
\\
A^{ijkm} & =\lambda
g^{ij}g^{km}+\mu(g^{ik}g^{jm}+g^{im}g^{jk}),\quad
\lambda=-\frac23\mu,%
\end{array}%
\right. \end{equation} where $g_{ij}$, and $g^{ij}$ are the
covariant and contravariant components of the metric tensor of
dimensional three Euclidian space in the curvilinear coordinate
$(x^1,x^2,\xi)$ define by \eqref{lab2.4} respectively,
\begin{equation}\label{lab3.3}
\begin{array}{ll}
\nabla_i w^j=\displaystyle\frac{\partial w^j}{\partial x^i}%
+\Gamma^j_{ik}w^k;\quad \nabla_i w_j=\displaystyle\frac{\partial w_j}{%
\partial x^i}-\Gamma^k_{ij}w_k, &  \\
\Gamma^i_{jk}=g^{il}(\displaystyle\frac{\partial g_{kl}}{\partial x^j}+%
\displaystyle\frac{\partial g_{jl}}{\partial x^k} -\displaystyle\frac{%
\partial g_{jk}}{\partial x^l}). &
\end{array}%
\end{equation}
 The absolute acceleration of the fluid is given by
\begin{equation}\label{lab3.4}
\begin{array}{ll}
a^i=\displaystyle\frac{\partial w^i}{\partial t}+w^j\nabla_j w^i
+2\varepsilon^{ijk}\omega_j w_k-\omega^2r^i; &  \\
a=\displaystyle\frac{\partial w}{\partial t}+(w\nabla)w+2\vec{\omega}\times%
\vec{w}+\vec{\omega}\times(\vec{\omega}\times\vec{R}), &
\end{array}%
\end{equation}
where $\vec{\omega}=\omega\vec{k}$ is the vector of angular velocity, $\vec{k%
}$ the unite vector along axis , and $\vec{R}$ the radium vector of
the fluid particle. The flow domain $\Omega_\varepsilon$ occupied by
the fluids in the channel in the impeller. The boundary $\partial
\Omega_\varepsilon$
of flow domain $\Omega_\varepsilon$ consists of inflow boundary $\Gamma_{in}$%
, out flow boundary $\Gamma_{out}$, positive blade's surface $S_+$,
negative blade's surface $S_-$ and top wall $\Gamma_t$ and Bottom
wall $\Gamma_b$:
\begin{equation}\label{lab3.5}
\partial\Omega_\varepsilon=\Gamma=\Gamma_{in}\cup\Gamma_{out}\cup S_-\cup
S_+\cup\Gamma_t\cup\Gamma_b \end{equation}
 Boundary conditions are given by
\begin{equation}\label{lab3.6}
\left\{%
\begin{array}{ll}
w|_{S_-\cup S_+}=0,\hskip 1cm w|_{\Gamma_b}=0,
w|_{\Gamma_t}=0,&\\
\sigma^{ij}(w,p)n_j|_{\Gamma_{in}}=g^i_{in},\quad
\sigma^{ij}(w,p)n_j|_{\Gamma_{out}}=g^i_{out} \text{(Natural
conditions)} &\\
\frac{\partial T}{\partial n}+\lambda(T-T_0)=0 \text{where
$\,\lambda\geq0$ is
constant,}%
\end{array}%
\right.\end{equation}
 If the fluid is incompressible and flow is stationary then
\begin{equation}\label{lab3.7}
\left\{%
\begin{array}{l}
\mbox{div} w=0, \\
(w\nabla)w+2\vec{\omega}\times\vec{w} +\nabla p-\nu\mbox{div}(e(w))=-\vec{%
\omega}\times(\vec{\omega}\times\vec{R}) +f, \\
w|_{\Gamma_0}=0,\quad \Gamma_0=S_+\cup S_-\cup \Gamma_t\cup \Gamma_b, \\
(-pn+2\nu e(w))|_{\Gamma_{in}}=g_{in},\quad \Gamma_1=\Gamma_{in}\cup
\Gamma_{0ut}, \\
(-pn+2\nu e(w))|_{\Gamma_{out}}=g_{out}, \\
w|_{t=0}=w_0(x),\quad \Omega_\varepsilon,%
\end{array}%
\right.
\end{equation}

For the polytropic ideal gas and flow is stationary, system
\eqref{lab3.1} turns to the conservation form
\begin{equation}\label{lab3.8}
\left\{%
\begin{array}{ll}
\mbox{div}(\rho w)=0, &  \\
\mbox{div}(\rho w\otimes w)+2\rho \omega\times w+R\nabla(\rho
T)&\\=\quad\mu \Delta w+(\lambda+\mu)\nabla\mbox{div} w-\rho
\omega\times(\omega\times \vec{R}), &
\\
\mbox{div}\lbrack \rho(\frac{|w|^2}{2}+c_vT+RT)w] &  \\
\quad=\kappa\Delta T+\lambda\mbox{div}(w\mbox{div}
w)+\mu\mbox{div}\lbrack w\nabla w]+\frac{\mu}{2}\Delta|w|^2, &
\end{array}%
\right.
\end{equation}
while for isentropic ideas gases, it turns
\begin{equation}\label{lab3.9}
\left\{%
\begin{array}{ll}
\mbox{div}(\rho w)=0, &  \\
\mbox{div}(\rho w\otimes w)+2\rho \omega\times
w+\alpha\nabla(\rho^\gamma)&\\
\quad=2\mu \mbox{div}(e) +\lambda\nabla\mbox{div}
w-\rho \omega\times(\omega\times \vec{R}), &  \\
&
\end{array}%
\right.
\end{equation}
where $\gamma>1 $is the specific heat radio and $\alpha$ a
positive constant.

The rate of work done by the impeller and global dissipative energy
are given by
\begin{equation}\label{lab3.10}
I(S,w(S))=\int\int_{S_-\cup S_+}\sigma\cdot n\cdot e_\theta\omega
r\mbox{d} S,\quad
J(S,w(S))=\int\int\int_{\Omega_\varepsilon}\Phi(w)\mbox{d} V
\end{equation}

where $e_\theta$ is base vector along the angular direction in a
polar cylindrical coordinate system.

Let us employ new coordinate system defined by \eqref{lab2.2}.
Flow's domain $\Omega_\varepsilon$ is mapped into $\Omega=D\times
[-1,1]$, where $D$ is a domain in $(x^1,x^2)\in \Re^2$ surround by
four arc $\widehat{AB},\widehat{CD},\widehat{CB},\widehat{DA}$ such
that
$$\partial D=\gamma_0\cup\gamma_1,\quad\gamma_0=\widehat{AB}\cup\widehat{CD},\quad
\gamma_1=\widehat{CB}\cup\widehat{DA},$$ and there exist four
positive functions
$\gamma_0(z),\widetilde{\gamma}_0(z),\gamma_1(z),\widetilde{\gamma}_1(z)$
such that
\begin{equation}\label{lab3.11}
\begin{array}{ll}
 r:=x^2=\gamma_0(x^1)=\gamma_0(z)\quad \mbox{on}\,\widehat{AB},\quad
   x^2=\widetilde{\gamma}_0(x^1)\quad \mbox{on}\,\widehat{CD}\\
 r:=x^2=\gamma_1(x^1)=\gamma_1(z)\,\quad \mbox{on}\,\widehat{DA},\quad
   x^2=\widetilde{\gamma}_1(x^1)\quad\mbox{on}\,\widehat{BC},\\
r_0\leq \gamma_0(z)\leq r_1\quad\mbox{on}\quad \widehat{AB},\quad
r_0\leq
\widetilde{\gamma}_0(z)\leq r_1\quad\mbox{on} \,\widehat{CD}, \\
 r_0\leq
\gamma_1(z)\leq r_1\quad\mbox{on}\, \widehat{DA},\quad r_0\leq
\widetilde{\gamma}_1(z)\leq r_1\quad\mbox{on}\,
\widehat{BC}.\end{array}
\end{equation}

%\begin{figure}[tbh]
%\begin{center}
%\includegraphics[bb=0 0 9cm 5cm]{MeridianPlane.bmp}
%\end{center}
%\caption{{\scriptsize Sectional Graph $D$ of Meridian Plane in
%Channel Flow $\Omega_\varepsilon$
%}}\label{MeridianPlane}
%\end{figure}

Let
\begin{equation}\label{lab3.12}
\begin{array}{ll}
\widetilde{\Gamma}_{in}=\vec{\Re}(\Gamma_{in}),\widetilde{\Gamma}_{out}=\vec{\Re}(\Gamma_{out}),
\widetilde{\Gamma}_{b}=\vec{\Re}(\Gamma_{b}),\widetilde{\Gamma}_{t}=\vec{\Re}(\Gamma_{t}),\\
\widetilde{\Gamma}_1=\widetilde{\Gamma}_{out}\cup\widetilde{\Gamma}_{in},\quad
\widetilde{\Gamma}_0=\widetilde{\Gamma}_{b}\cup\widetilde{\Gamma}_{t}\cup\{\xi=1\}\cup\{\xi=-1\},
\end{array}
\end{equation}
\begin{equation}\label{lab3.13}
\begin{array}{ll}
\partial D=\gamma_0\cup\gamma_1, \quad\partial
\Omega=\widetilde{\Gamma}_0\cup\widetilde{\Gamma}_1,\\
\gamma_0=(D\cap\widetilde{\Gamma}_b)\cup(D\cap\widetilde{\Gamma}_t),
\quad\gamma_1=(D\cup\widetilde{\Gamma}_{out})\cup(D\cup\widetilde{\Gamma}_{in}),
\end{array}
\end{equation}
 where $\vec{\Re}$ is defined by \eqref{lab2.1}.

 Let denote
\begin{equation}\label{lab3.14}
V(\Omega):=\{v|,v\,\in%
\,H^1(\Omega)^3,\,v|_{\widetilde{\Gamma}_0}=0\},H^1_\Gamma(\Omega)=
\{q|,q\,\in\,H^1(\Omega),\,q|_{\widetilde{\Gamma}_0}=0\}.
\end{equation}
 The
variational formulations for Navier-Stokes problem \eqref{lab3.7}
and \eqref{lab3.9} are respectively given by
\begin{equation}\label{lab3.15}
\left\{%
\begin{array}{ll}
\mbox{Find}\,(w,p),w\,\in\,V(\Omega),p\in\,L^2(\Omega),\mbox{such
that} &
\\
a(w,v)+2(\omega\times w,v)+b(w,w,v)+ &  \\
\quad\quad -(p, \mbox{div} v)=\langle F,v\rangle,\quad \forall\,v\,\in\,V(\Omega), &  \\
(q,\mbox{div} w)=0,\quad \forall\,q\,\in\,L^2(\Omega), &  \\
&
\end{array}%
\right.
\end{equation}
and
\begin{equation}\label{lab3.16}
\left\{%
\begin{array}{ll}
\mbox{Find}\,(w,\rho),w\,\in\,V(\Omega),\rho\in\,L^\gamma(\Omega),\mbox{such
that} &  \\
a(w,v)+2(\omega\times w,v)+b(\rho w,w,v)+ &  \\
\quad\quad +(-\alpha p+\lambda\mbox{div} w, \mbox{div} v)=\langle
F,v\rangle,\quad
\forall\,v\,\in\,V(\Omega), &  \\
(\nabla q,\rho w)=\langle \rho wn,q\rangle|_{\Gamma_1},\quad
\forall\,q\,\in\,H^1_\Gamma(\Omega), &  \\
&
\end{array}%
\right.
\end{equation}
 where
\begin{equation}\label{lab3.17}
\begin{array}{ll}
\langle F,v\rangle:=<f,v>+<\widetilde{g},v>_{\widetilde{\Gamma}_{1}},&\\
\langle \widetilde{g}%
,v\rangle=\langle g_{in},v\rangle|_{\widetilde{\Gamma}_{in}}+\langle g_{out},v\rangle|_{\widetilde{\Gamma}_{out}}, &  \\
a(w,v)=\int_{\Omega}A^{ijkm}e_{ij}(w)e_{km}(v)\sqrt{g}dxd\xi, &  \\
b(w,w,v)=\int_\Omega g_{km}w^j\nabla_jw^kv^m\sqrt{g}dxd\xi, &
\end{array}%
\end{equation}
Next we rewrite \eqref{lab3.7} and \eqref{lab3.9} in new coordinate
system. Because second kind of Christoffel symbols in new coordinate
system are
\begin{equation}\label{lab3.18}
\left\{%
\begin{array}{l}
\Gamma^\alpha_{\beta\gamma}=-r\delta_{2\alpha}\Theta_\beta\Theta_\gamma,
\quad\Gamma^\alpha_{3\beta}= -\varepsilon r\delta_{2\alpha}\Theta_\beta, \\
\Gamma^3_{\alpha\beta}=\varepsilon^{-1}r^{-1}(\delta_{2\alpha}\delta^%
\lambda_\beta+ \delta_{2\beta}\delta^\lambda_\alpha
)\Theta_\lambda+\varepsilon^{-1}\Theta_{\alpha\beta}+\varepsilon^{-1}r%
\Theta_2\Theta_\alpha\Theta_\beta, \\
\Gamma^3_{3\alpha}=\Gamma^3_{\alpha
3}=r^{-1}\delta_{2\alpha}+r\Theta_2\Theta_\alpha\quad
\Gamma^\alpha_{33}=-\varepsilon^2
r\delta_{2\alpha},\quad\Gamma^3_{33}=\varepsilon r\Theta_2.%
\end{array}%
\right.
\end{equation}
 the covariant derivatives of the velocity field
$\nabla_iw^j=\frac{\partial w^j}{\partial x^i}+\Gamma^j_{ik}w^k$ can
be expressed as
\begin{lem}\label{lem1}
 Under the curvilinear coordinate system  $%
(x^1,x^2,\xi)$ defined by \eqref{lab2.4}, the covariant derivatives
of the velocity field can  be expressed as
\begin{equation}\label{lab3.19}
\left\{%
\begin{array}{l}
\nabla_\alpha w^\beta=\frac{\partial w^\beta}{\partial x^\alpha}
-r\delta^\beta_{2}\Theta_\alpha\Pi(w,\Theta), \\
\nabla_\alpha w^3=\frac{\partial w^3}{\partial x^\alpha}%
+\varepsilon^{-1}(x^2)^{-1}w^2\Theta_\alpha+
\varepsilon^{-1}w^\beta\Theta_{\alpha\beta}+(\varepsilon
x^2)^{-1}a_{2\alpha}\Pi(w,\Theta), \\
\nabla_3 w^\alpha=\frac{\partial
w^\alpha}{\partial\xi}-x^2\varepsilon
\delta_{2\alpha}\Pi(w,\Theta), \\
\nabla_3 w^3=\frac{\partial w^3}{\partial \xi}+\frac{w^2}{x^2}%
+x^2\Theta_2\Pi(w,\Theta), \\
\mbox{div} w=\frac{\partial w^\alpha}{\partial x^\alpha}+\frac{w^2}{x^2}+%
\frac{\partial w^3}{\partial\xi},\quad \Pi(w,\Theta)=\varepsilon
w^3+w^\beta\Theta_\beta.%
\end{array}%
\right.\end{equation}
 and the deformation tensors are given by
\begin{equation}\label{lab3.20}
\left\{
\begin{array}{ll}
e_{ij}(w)=\varphi_{ij}(w)+\psi_{ij}(w,\Theta);\\
\psi_{ij}(w,\Theta)=\psi_{ij}^\lambda(w)\Theta_\lambda+\psi_{ij}^{\lambda\sigma}(w)\Theta_\lambda\Theta_\sigma
+e_{ij}^*(w,\Theta),
\end{array}\right.
\end{equation}
where
\begin{equation}\label{lab3.21}
\left\{%
\begin{array}{rl}
\varphi_{\alpha\beta}(w) &
=\displaystyle\frac12(\displaystyle\frac{\partial w^\alpha}{\partial
x^\beta} +\displaystyle\frac{\partial w^\beta}{\partial
x^\alpha}),\\
\varphi_{3\alpha}(w) &
=\displaystyle\frac12(\displaystyle\frac{\partial
w^\alpha}{\partial\xi}+\varepsilon^{2}r^2\displaystyle\frac{\partial w^3}{%
\partial x^\alpha}),\quad
\varphi_{33}(w)  =\varepsilon^2r^2(\displaystyle\frac{\partial w^3}{%
\partial\xi}+\displaystyle\frac{w^2}{r}),\\
\end{array}\right.
\end{equation}
\begin{equation}\label{lab3.22}
\left\{\begin{array}{ll}
\psi^\lambda_{\alpha\beta}(w)&=\displaystyle%
\frac12\varepsilon r^2(\displaystyle\frac{\partial w^3}{\partial x^\alpha}%
\delta^\lambda_\beta+\displaystyle\frac{\partial w^3} {\partial x^\beta}%
\delta^\lambda_\alpha), \\
 \psi^\lambda_{3\alpha}(w)&=\displaystyle\frac12
\varepsilon r^2(\displaystyle\frac{\partial w^\lambda}{\partial x^\alpha}%
+\delta^\lambda_\alpha(\displaystyle\frac{\partial w^3}{\partial\xi} +%
\displaystyle\frac{2}{r}w^2)), \quad
\psi^\lambda_{33}(w)=\varepsilon r^2\frac{\partial
w^\lambda}{\partial\xi},\\
\psi^{\lambda\sigma}_{\alpha\beta}(w) & =\displaystyle\frac12r^2 (%
\displaystyle\frac{\partial w^\lambda}{\partial x^\alpha}\delta_{\beta%
\sigma} +\displaystyle\frac{\partial w^\lambda}{\partial x^\beta}%
\delta_{\sigma\alpha}+\frac2{r}w^2\delta_{\alpha\lambda}\delta_{\sigma%
\beta}) , \\
\psi^{\lambda\sigma}_{3\alpha}(w) & =\frac12r^2\displaystyle\frac{%
\partial w^\lambda}{\partial \xi}\delta_{\alpha\sigma}, \quad
\psi^{\lambda\sigma}_{33}(w)=0.%
\end{array}%
\right.
\end{equation}
\begin{equation}\label{lab3.23}
\begin{array}{ll}
e^*_{\alpha\beta}(w,\Theta)=\frac12r^2w^\sigma(\Theta_\alpha\Theta_{\sigma%
\beta}+\Theta_\beta\Theta_{\sigma\alpha}),
e^*_{3\alpha}(w)=\frac12\varepsilon
r^2w^\sigma\Theta_{\sigma\alpha}, e^*_{33}(w)=0.
\end{array}
\end{equation}
\end{lem}
The proof is omitted here.

 By simply tensor calculations in terms of
\eqref{lab3.19} \eqref{lab2.5} and \eqref{lab2.6},
\begin{equation}\label{lab3.24}
\begin{array}{ll}
 &A^{ijkl}e_{kl}(w)e_{ij}(v)\\&=e_{\alpha\beta}(w)e_{\alpha\beta}(v)
 +2(g^{33}\delta^{\alpha\beta}+\varepsilon^{-2}\Theta_\alpha\Theta_%
\beta)e_{3\alpha}(w)e_{3\beta}(v) \\
& +g^{33}g^{33}e_{33}(w)e_{33}(v)
-2\varepsilon^{-1}\Theta_\beta(e_{\alpha\beta}(w)e_{3\alpha}(v)
+e_{\alpha\beta}(v)e_{3\alpha}(w)) \\
& -2\varepsilon^{-1}
g^{33}\Theta_\alpha(e_{33}(w)e_{3\alpha}(v)+e_{33}(v)e_{3\alpha}(w)) \\
& +\varepsilon^{-2}\Theta_\alpha\Theta_\beta(e_{\alpha\beta}(w)e_{33}(v)+e_{%
\alpha\beta}(v)e_{33}(w))
=A^{ijkl}\varphi_{kl}(w)\varphi_{ij}(v)\\
&+A^{ijkl}[\varphi_{kl}(w)\psi_{ij}(v,\Theta)+\psi_{kl}(w,\Theta)\varphi_{ij}(v)
+\psi_{kl}(w,\Theta)\psi_{ij}(v,\Theta)],
\end{array}\end{equation}
 and
$$\int_{\Omega}A^{ijkl}\varphi_{kl}(w)\varphi_{ij}(v)r\varepsilon d\xi dx=((w,v))+\int_{\Omega}B(w,v,\Theta)d\xi dx,$$
where
\begin{equation}\label{lab3.25}
\left\{\begin{array}{ll}
((w,v))=\int\limits_\Omega2\mu\left[\varphi_{\alpha\beta}(w)\varphi_{
\alpha\beta}(v)
+2(r\varepsilon)^{-2}\varphi_{3\alpha}(w)\varphi_{3\beta}(v)
+(r\varepsilon)^{-4}\varphi_{33}(w)\varphi_{33}(v)\right]r\varepsilon
d\xi
dx,\\
\|w\|^2_\Omega:=((w,w)),\\
a(w,v)
=((w,v))+\int\limits_\Omega2\nu\left[B(w,v,\Theta)+A(w,v,\Theta)\right]r\varepsilon
 d\xi dx,
\end{array}\right.
\end{equation}
and
\begin{equation}\label{lab3.26}
\left\{\begin{array}{ll}
 B(w,v,\Theta)&=2\mu[2\varepsilon^{-2}(|\nabla\Theta|^2\delta_{\alpha\beta}+\Theta_\alpha\Theta_\beta)\varphi_{3\alpha}(w)
 \varphi_{3\beta}(v)\\
 &+\varepsilon^{-4}r^{-2}|\nabla\Theta|^2(1+a)\varphi_{33}(w)\varphi_{33}(v)\\
& -2\varepsilon^{-1}\Theta_\beta(\varphi_{\alpha\beta}(w)\varphi_{3\alpha}(v)+\varphi_{\alpha\beta}(v)\varphi_{3\alpha}(w))\\
&-2\varepsilon^{-1}g^{33}\Theta_\alpha(\varphi_{33}(w)\varphi_{3\alpha}(v)+\varphi_{33}(v)\varphi_{3\alpha}(w))\\
&+\varepsilon^{-2}\Theta_\alpha\Theta_\beta(\varphi_{\alpha\beta}(w)\varphi_{33}(v)
+\varphi_{\alpha\beta}(v)\varphi_{33}(w))],\\
A(w,v,\Theta):&=A^{ijkl}[\varphi_{kl}(w)\psi_{ij}(v,\Theta)
+\varphi_{ij}(v)\psi_{kl}(w,\Theta)+\psi_{ij}(w,\Theta)\psi_{ij}(v,\Theta)]
\end{array}\right.
\end{equation}
In particular,
\begin{equation}\label{lab3.27}
\left\{\begin{array}{ll}
B(w,w,\Theta)&=2\mu[2\varepsilon^{-2}(|\nabla\Theta|^2\delta_{\alpha\beta}+\Theta_\alpha\Theta_\beta)
\varphi_{3\alpha}(w) \varphi_{3\beta}(w)\\
 &+\varepsilon^{-4}r^{-2}|\nabla\Theta|^2(1+a)\varphi_{33}(w)\varphi_{33}(w)
 -4\varepsilon^{-1}\Theta_\beta\varphi_{\alpha\beta}(w)\varphi_{3\alpha}(w)\\
&-4\varepsilon^{-1}g^{33}\Theta_\alpha\varphi_{33}(w)\varphi_{3\alpha}(w)
+2\varepsilon^{-2}\Theta_\alpha\Theta_\beta\varphi_{\alpha\beta}(w)\varphi_{33}(w)],\\
 A(w,w,\Theta)&=A^{ijkl}[2\varphi_{kl}(w)\psi_{ij}(w,\Theta)
+\psi_{ij}(w,\Theta)\psi_{ij}(w,\Theta)]
\end{array}\right.
\end{equation}

Next we consider trilinear form. Using \eqref{lab3.19}, we have
\[
\begin{array}{ll}
w^j\nabla_jw^\beta=w^\alpha\frac{\partial w^\beta}{\partial x^\alpha}+w^3%
\frac{\partial w^\beta}{\partial \xi}-r\delta_{2\beta}\Pi(w,\Theta)\Pi(w,%
\Theta), &  \\
w^j\nabla_jw^3=w^\alpha\frac{\partial w^3}{\partial x^\alpha}+w^3(\frac{%
\partial w^3}{\partial\xi}+\frac{w^2}{r})+\varepsilon^{-1}w^\alpha
w^\beta\Theta_{\alpha\beta} &  \\
\quad+\varepsilon^{-1}r\Theta_2\Pi(w,\Theta)\Pi(w,\Theta)
+\varepsilon^{-1}r^{-1}w^2(w^\alpha\Theta_\alpha+\Pi(w,\Theta)), &  \\
&
\end{array}%
\]
Therefore
\begin{equation}\label{lab3.28}
\begin{array}{ll}
&b(w,w,v) \\& =\int_D\int^1_{-1}g_{km}w^j\nabla_jw^kv^m\sqrt{g}d\xi
dx
=\int_D\int^1_{-1}\{[a_{\alpha\beta}(w^\lambda\frac{\partial w^\alpha}{%
\partial x^\lambda}+w^3\frac{\partial w^\alpha}{\partial\xi}) \\
& +\varepsilon r^2\Theta_\beta(w^\lambda\frac{\partial w^3}{\partial
x^\lambda}+w^3(\frac{\partial w^3}{\partial\xi}+\frac{w^2}{r}%
))-r\delta_{2\beta}\Pi(w,\Theta)\Pi(w,\Theta) \\
& +r^2w^\sigma w^\lambda\Theta_{\sigma\lambda}\Theta_\beta
+r\Theta_\beta
w^2(w^\lambda\Theta_\lambda+\Pi(w,\Theta))]v^\beta \\
& + \varepsilon r^2[\Theta_\alpha(w^\lambda\frac{\partial
w^\alpha}{\partial
x^\lambda}+w^3\frac{\partial w^\alpha}{\partial\xi})+\varepsilon(w^\lambda%
\frac{\partial w^3}{\partial x^\lambda}+w^3 (\frac{\partial w^3}{\partial\xi}%
+\frac{w^2}{r})) \\
& +w^\alpha w^\beta\Theta_{\alpha\beta}
+w^2(w^\alpha\Theta_\alpha+\Pi(w,\Theta))]v^3\}r\varepsilon d\xi dx%
\end{array}%
\end{equation}
 By similar manner,
angular velocity vector is given by
\begin{equation}\label{lab3.29}
\left\{%
\begin{array}{ll}
\vec{\omega}=\omega \vec{e_1}-\varepsilon^{-1}\omega\Theta_1\vec{e_3}, &  \\
2\vec{\omega}\times \vec{w}&\\=-2r
\omega\delta_{\alpha2}(\varepsilon
w^3+w^\beta\Theta_\beta)\vec{e_\alpha} +2r\omega\varepsilon^{-1}
[\Theta_2(\varepsilon w^3+w^\beta\Theta_\beta)+r^{-2}w^2]\vec{e_3} &  \\
=-2r \omega\delta_{\alpha2}\Pi(w,\Theta)\vec{e_\alpha}
+2r\omega\varepsilon^{-1}
[\Theta_2\Pi(w,\Theta)+r^{-2}w^2]\vec{e_3}, &
\end{array}%
\right.
\end{equation}
 while Coliali form
\begin{equation}\label{lab3.30}
\begin{array}{ll} &C(w,v) \\&
:=\int_D\int^1_{-1}2g_{ij}(\vec{\omega}\times w)^iv^j\sqrt{g}d\xi dx
\\
& =\int_D\int^1_{-1}[2a_{\alpha\beta}(\vec{\omega}\times w)^\alpha
v^\beta+2\varepsilon(x^2)^2\Theta_\alpha((\vec{\omega}\times w)^\alpha v^3+(%
\vec{\omega}\times w)^3v^\alpha) \\
& +2\varepsilon^2(x^2)^2(\vec{\omega}\times w)^3v^3]\varepsilon
x^2d\xi dx
\\
& =\int_D\int^1_{-1}2r\omega[(w^2\Theta_\beta-\delta_{2\beta}\Pi(w,%
\Theta))v^\beta+\varepsilon w^2v^3]r\varepsilon d\xi dx%
\end{array}
\end{equation}

It is very easy to verify that
\begin{equation}\label{lab3.31}
C(w,w)=0.
\end{equation}

Throughout this paper, Latin indices and exponents $i,j,k \cdots$
vary in the set $\{1,2,3\}$, while Greek indices and exponents
$\alpha,\beta,\gamma\cdots$ vary in the set $\{1,2\}$ . Furthermore,
the summation convention with respect to repeated indices or
exponents is systematically used in conjunction with this rule.
\begin{lem}\label{lem2}
The function $\|\cdot\|_\Omega$ defined by (3.25) is a norm in
Hilbert space
$$V(\Omega):=\{v\in\,H^1(D)^3,v|_{\widetilde{\Gamma}_1}=0,\}\eqno{(3.32)}$$
\end{lem}
\begin{proof}
Indeed it is enough to prove that
$$\|w\|_\Omega=0 ,\forall\,w\in\,V(\Omega)\Rightarrow w=0.$$
This means that
$$\|w\|_\Omega=0,\mbox{i.e., } \varphi_{ij}(w)=0.$$ we have to prove
$w=0$. Firstly, the following identity is held
$$\partial_\gamma(\partial_\alpha
w^\beta)=\partial_\gamma\varphi_{\alpha\beta}(w)+\partial_\alpha\varphi_{\gamma\beta}(w)
-\partial_\beta\varphi_{\alpha\gamma}(w).$$ This shows that
$$\varphi_{\alpha\beta}(w)=0,\, \mbox{in}\,D\Rightarrow
\partial_\gamma\partial_\alpha w^\beta=0,\,\mbox{in}\,{\cal D}'(D).$$
By a classical result from distribution theory, each function $w$ is
therefore a polynomial of degree$\leq1$(recall that the set $D$ is
connect). In other words, there exist constants $c_\alpha$ and
$d_{\alpha\beta}$ such that
$$w^\alpha(x)=c_\alpha+d_{\alpha\beta}x^\beta,\,\forall\,x=(x^1,x^2)\in\,D,$$
But $\varphi_{\alpha\beta}(w)=0$ also implies that
$d_{\alpha\beta}=-d_{\beta\alpha}$; hence there exist two vectors
$\vec{c},\vec{d}\in\,\Re^2$ such that
$$w=\vec{c}+\vec{d}\times \vec{x},\forall\,x\,\in\,D,$$
Since $w|_{\widetilde{\Gamma}}=0$ and the set where such a vector
field $w^\alpha$ vanishes is always of zero area unless
$\vec{c}=\vec{d}=0$, it follows that $w^\alpha=0$ when area
$\widetilde{\Gamma}_0>0$. On the other hand, in view of boundary
condition (3.13)
$$\varphi_{33}(w)  =\varepsilon^2r^2(\displaystyle\frac{\partial w^3}{%
\partial\xi}+\displaystyle\frac{w^2}{r})=0\Rightarrow \frac{\partial w^3}{%
\partial\xi}=0\Rightarrow w^3=0,$$
The proof is complete.
\end{proof}

\begin{lem}\label{lem3}
The norm $\|\cdot\|_\Omega $ and the norm
$$|w|_{1,\Omega}^2=\int\limits_\Omega[\sum\limits_{i=1}^3(\sum\limits_{\alpha=1}^2(\frac{\partial
w^i}{\partial x^\alpha})^2+(\frac{\partial
w^i}{\partial\xi})^2)]r\varepsilon d\xi dx, \quad\forall
w\,\in\,V(\Omega),$$ are equivalent in $V(\Omega)$,i.e. there exist
a constant $C_i(\Omega)>0,i=1,2$ depending upon $\Omega$ only such
that
$$ C_1(\Omega)|w|_{1,\Omega}\leq \|w\|_\Omega\leq
C_2(\Omega)|w|_{1,\Omega},\quad\forall\,w\,\in\,V(\Omega),\eqno{(3.33)}$$
\end{lem}

\begin{proof} Firstly we indicate that in view of (3.11)(3.12) we assert
that there exist a constant $C_i(\Omega)>0,i=1,2$ depending upon
$\Omega$ only such that
$$ C_1(\Omega)(\sum\limits_{i,j=1}^3\|\varphi_{ij}(w)\|^2_{0,\Omega})^{1/2}\leq \|w\|_\Omega\leq
C_2(\Omega)(\sum\limits_{i,j=1}^3\|\varphi_{ij}(w)\|^2_{0,\Omega})^{1/2},\quad\forall\,w\,\in\,V(\Omega),$$
and $\varphi_{ij}(w)$ can be looked  as strain tensor in Cartesian
coordinates in $\Re^3$ then according to Korn's inequality (see
[14][15]) the
$(\sum\limits_{i,j=1}^3\|\varphi_{ij}(w)\|^2_{0,\Omega})^{1/2}$ is
equivalent to $\|w\|_{1,\Omega}$, therefor this reach to (3.33). The
proof is complete.
\end{proof}

\begin{thm}\label{thm3}
 Under the new coordinate system, the stationary
Navier-Stokes equations can be explicitly expressed as $\Theta$:

$$\left\{\begin{array}{ll}
\mbox{div}(w)=\frac{\partial w^\alpha}{\partial x^\alpha}+\frac{\partial w^3%
}{\partial \xi}  +\frac{w^2}{r}=0,\\

\mathcal{N}^i(w,p,\Theta):=\mathcal{L}^i(w,p,\Theta)+\mathcal{NN}%
^i(u,u)=f^i,\quad i=1,2,3,  \\

\end{array}\right.\eqno{(3.34)}$$

where
$$
\begin{array}{ll}
\mathcal{L}^\alpha(w,p,\Theta):=-\nu[\widetilde{\Delta} w^\alpha-2%
\varepsilon^{-1}\Theta_\beta\frac{\partial^2w^\alpha}{\partial\xi\partial
x^\beta}+(r\varepsilon)^{-2}a\frac{\partial^2w^\alpha}{\partial(\xi)^2}
+r^{-1}\frac{\partial w^\alpha}{\partial x^2} &  \\
\quad-2r\varepsilon\delta_{2\alpha}\Theta_\lambda\frac{\partial w^3}{%
\partial x^\lambda} -(r\varepsilon)^{-1}(\delta_{\alpha\lambda}\Theta_2+2%
\delta_{2\alpha}\Theta_\lambda +r\delta_{\alpha\lambda}\widetilde{\Delta}%
\Theta)\frac{\partial w^\lambda}{\partial\xi}-2r^{-1}\delta_{2\alpha}\frac{%
\partial w^3}{\partial\xi} &  \\
\quad-2\varepsilon\delta_{2\alpha}a\Theta_2w^3+
\delta_{2\alpha}(r^{-2}\delta_{2\sigma}-2\Theta_2\Theta_\sigma -a_{2\sigma}|%
\widetilde{\nabla}\Theta|^2-r\Theta_\lambda\Theta_{\lambda\sigma})w^\sigma)
]
&  \\
\quad-2r\delta_{2\alpha}\omega\Pi(w,\Theta)+\nabla_\alpha p
-\varepsilon^{-1}\Theta_\alpha\frac{\partial p}{\partial\xi}, &  \\
\mathcal{NN}^\alpha(w,w):=w^\beta\frac{\partial w^\alpha}{\partial x^\beta}%
+w^3\frac{\partial w^\alpha}{\partial\xi}-r\delta_{2\alpha}\Pi(w,\Theta)%
\Pi(w,\Theta), &  \\
&
\end{array}%
\eqno{(3.35)}
$$
$$
\begin{array}{ll}
\mathcal{L}^3(w,p,\Theta):=-\nu[\widetilde{\Delta} w^3-2\varepsilon^{-1}%
\Theta_\beta\frac{\partial^2 w^3}{\partial\xi\partial x^\beta}%
+(r\varepsilon)^{-2}a\frac{\partial^2 w^3}{\partial\xi^2} +3r^{-1}\delta_{2%
\beta}\frac{\partial w^3}{\partial x^\beta} &  \\
\quad +\varepsilon^{-2}r^{-3}a \frac{\partial w^2}{\partial\xi}%
+2(r\varepsilon)^{-1}(\delta_{2\beta}\Theta_\lambda +r\Theta_{\beta\lambda})%
\frac{\partial w^\lambda}{\partial x^\beta} &  \\
\quad-(r\varepsilon)^{-1}(\Theta_2+r\widetilde{\Delta}\Theta)\frac{\partial
w^3}{\partial \xi}+(r\varepsilon)^{-1}(r^{-1}\delta_{2\alpha}\Theta_2+3%
\Theta_{2\alpha} +r\partial_\alpha\widetilde{\Delta}\Theta)w^\alpha] &  \\
\quad-2r\omega\varepsilon^{-1}\Pi(w,\Theta)\Theta_2 -2(\varepsilon
r)^{-1}w^2 -\varepsilon^{-1}\Theta_\beta\partial_\beta p+(r\varepsilon)^{-1}a%
\frac{\partial p}{\partial\xi}, &  \\
\mathcal{NN}^\alpha(w,w):=w^\beta\frac{\partial w^3}{\partial x^\beta}+w^3%
\frac{\partial w^3}{\partial\xi}+\varepsilon^{-1}w^\beta
w^\lambda\Theta_{\beta\lambda}+(r\varepsilon)^{-1}\Pi(w,\Theta)(2w^2 &  \\
\quad+r^2\Theta_2\Pi(w,\Theta)), &
\end{array}%
\eqno{(3.36)}
$$
where $\widetilde{\Delta}\varphi=\frac{\partial^2\varphi}{\partial
x^1\partial x^1}+\frac{\partial^2\varphi}{\partial x^2\partial x^2}$,\quad $%
\hat{\varepsilon}^{\alpha\beta}=+1,-1,0 $ depending on
$(\alpha,\beta)= (1,2),(2,1),\mbox{or other case}$. Under the new
coordinate system the domain $\Omega=\{(x^1,x^2)\in\,D,-1\leq
\xi\leq 1\}$ is independent  of geometry of blade's surface
$\Theta$.
\end{thm}

\begin{proof} Substituting (3.19)(2.5)(2.6) and (3.20) into (3.7), tensor
calculations show that (3.14) is valid. The details is omitted here.
\end{proof}

By using above formulae we claim that there exist of the Gateaux
derivative of the solution of Navier-Stokes equations to be satisfy
following linearized Navier-Stokes equation

\begin{thm}\label{thm4}Assume that there exists a solution $%
(w(\Theta),p(\Theta))$  of Navier-Stokes problem (3.7) such that
define a mapping $\Theta\Rightarrow (w(\Theta),p(\Theta))$ from
$H^1_0(D)\cap H^2(D)$
to $H^{1,q}(\Omega)\times L^{2,q}(\Omega)$. Then  there exist the G$\hat{a}$%
teaux derivatives of $(w,p)$ at a point $\Theta\,\in H^1_0(D)\cap
H^2(D)$
with respect to any direction $\eta\in H^1_0(D)\cap H^2(D)$:$\hat{w}\doteq%
\displaystyle\frac{\mathcal{D}w}{\mathcal{D}\Theta}\eta,  \hat{p}\doteq%
\displaystyle\frac{\mathcal{D}p}{\mathcal{D}\Theta}\eta $ and
satisfy following linearized equations:
$$
\left\{%
\begin{array}{l}
\mbox{div}(\hat{w})=0 \\
\mathcal{L}^i(\hat{w},\widehat{p},\Theta)+\mathcal{NN}^i(w,\widehat{w})+%
\mathcal{NN}^i(\widehat{w},w)=R^i(w,p) \\
\hat{w}|_{\Gamma_s}=0,\quad \sigma^{ij}(\hat{w},\hat{p})n_j|_{\Gamma_0}=0.%
\end{array}%
\right.\eqno(3.37)
$$%
where
$$
\begin{array}{ll}
R^\alpha(w,p,\Theta) &
=\nu\delta_{2\alpha}\frac{\partial^2}{\partial
x^\lambda\partial x^\sigma} (rw^\sigma\Theta_\lambda)+\nu\varepsilon^{-1}%
\widetilde{\Delta}\frac{\partial w^\alpha}{\partial\xi} \\
& +\nu\frac{\partial}{\partial x^\beta} [-2(\varepsilon)^{-2}\Theta_\beta%
\frac{\partial^2w^\alpha}{\partial\xi^2} +2r\varepsilon\delta_{2\alpha}\frac{%
\partial w^3}{\partial x^\beta} \\
& +(r\varepsilon)^{-1}(\delta_{\alpha\lambda}\delta_{2\beta}+2\delta_{2%
\alpha}\delta_{\lambda\beta}) \frac{\partial w^\lambda}{\partial\xi}
+2\varepsilon\delta_{2\alpha}w^3(a\delta_{2\beta}+2r^2\Theta_\beta\Theta_2)
\\
& +
\delta_{2\alpha}w^\sigma(2\delta_{2\beta}\Theta_\sigma+2\Theta_2\delta_{%
\beta\sigma}
+2r^2(\Theta_2\delta_{\sigma\beta}+\delta_{2\beta}\Theta_\sigma)|\widetilde{%
\nabla}\Theta|^2 \\
& -2a_{2\sigma}\Theta_\beta -r\Theta_{\beta\sigma})] -\frac{\partial}{%
\partial x^\beta}[2r\delta_{2\alpha}(\Pi(w,\Theta) +\omega)w^\beta
+\varepsilon^{-1}\delta_{\alpha\beta}\frac{\partial p}{\partial\xi}] \\
&
\end{array}%
\eqno{(3.38)}
$$
$$
\begin{array}{ll}
R^3(w,p,\Theta) & =\nu\varepsilon^{-1}(-\widetilde{\Delta}(\frac{w^2}{r}%
)+3\partial_{2\alpha}(r^{-1}w^\alpha))
+\varepsilon^{-1}\partial_{\beta\lambda}(w^\beta w^\lambda) -2\nu\frac{%
\partial}{\partial x^\beta}(\Theta_\beta\frac{\partial^2 w^3}{\partial \xi^2}%
) \\
& -\nu(\varepsilon)^{-1}\frac{\partial}{\partial x^\beta}[2(r%
\varepsilon)^{-1}\Theta_\beta\frac{\partial w^2}{\partial \xi}+2r^{-1}\frac{%
\partial w^\beta}{\partial x^2}-r^{-1}\delta_{2\beta}\frac{\partial w^3}{%
\partial \xi}+r^{-2}\delta_{2\beta}w^2] \\
& +\varepsilon^{-1}\frac{\partial}{\partial x^\beta}[2w^\beta\frac{w^2}{r}%
-2\omega r\omega\delta_{2\beta}\Pi(w,\Theta) \\
&
+r(\Pi(w,\Theta)-\omega)(2\Theta_2w^\beta-\delta_{2\beta}\Pi(w,\Theta))]
+\varepsilon^{-1}\widetilde{\Delta}p+2\varepsilon^{-1}\frac{\partial}{%
\partial x^\beta}(r\Theta_\beta\frac{\partial p}{\partial \xi}) \\
&
\end{array}%
\eqno{(3.39)}
$$
\end{thm}

\begin{proof}The Navier-Stokes equations (3.24) read

$$
\begin{array}{ll}
\frac{\partial w^\alpha}{\partial x^\alpha}+\frac{w^2}{r}+\frac{\partial w^3%
}{\partial x^3}=0, &  \\
\mathcal{N}^\alpha(w,p,\Theta) \vec{e}_\alpha+\mathcal{N}^3(w,p,\Theta)\vec{e%
}_3=f^\alpha \vec{e}_\alpha+f^3 \vec{e}_3. &
\end{array}%
\eqno{(3.40)}
$$
Set Gateaux derivative with respect with $\Theta$ along any director $%
\eta\in\,\mathcal{W}:=H^2(D)\cap H^1_0(D)$ denoted by $\frac{\mathcal{D}}{%
\mathcal{D}\Theta}\eta$. Then from (3.40) we obtain
\[
\begin{array}{ll}
\frac{\mathcal{D}}{\mathcal{D}\Theta}\mathcal{N}^\alpha(w,p,\Theta) \vec{e}%
_\alpha\eta
+\frac{\mathcal{D}}{\mathcal{D}\Theta}\mathcal{N}^3(w,p,\Theta)
\vec{e}_3\eta+\mathcal{N}^\alpha(w,p,\Theta)\frac{\mathcal{D}\vec{e}_\alpha}{%
\mathcal{D}\Theta} \eta+\mathcal{N}^3(w,p,\Theta)\frac{\mathcal{D}\vec{e}_3}{%
\mathcal{D}\Theta}\eta &  \\
\quad=f^\alpha \frac{\mathcal{D}\vec{e}_\alpha}{\mathcal{D}\Theta}\vec{e}%
_\alpha\eta+f^3
\frac{\mathcal{D}\vec{e}_3}{\mathcal{D}\Theta}\vec{e}_3\eta,
&  \\
\frac{\mathcal{D}}{\mathcal{D}\Theta}\mathcal{N}^\alpha(w,p,\Theta) \vec{e}%
_\alpha +\frac{\mathcal{D}}{\mathcal{D}\Theta}\mathcal{N}^3(w,p,\Theta) \vec{%
e}_3+[\mathcal{N}^\alpha(w,p,\Theta)-f^\alpha]\frac{\mathcal{D}\vec{e}_\alpha%
}{\mathcal{D}\Theta}+[\mathcal{N}^3(w,p,\Theta)-f^3]\frac{\mathcal{D}\vec{e}%
_3}{\mathcal{D}\Theta}\vec{e}_3=0, &  \\
&
\end{array}%
\]
Hence,
\[
\begin{array}{ll}
\frac{\mathcal{D}}{\mathcal{D}\Theta}\mathcal{N}^\alpha(w,p,\Theta)\eta=%
\mathcal{L}^\alpha(\widehat{w},\widehat{p},\Theta)\eta+\mathcal{NN}^\alpha(w,%
\widehat{w})\eta+\mathcal{NN}^\alpha(\widehat{w},w)\eta+R^\alpha(w,p,\Theta)%
\eta=0, &  \\
\frac{\mathcal{D}}{\mathcal{D}\Theta}\mathcal{N}^3(w,p,\Theta)\eta=\mathcal{L%
}^3(\widehat{w},p,\Theta)\eta++\mathcal{NN}^3(w,\widehat{w})\eta+\mathcal{NN}%
^3(\widehat{w},w)\eta+R^3(w,p,\Theta)\eta=0, &  \\
&
\end{array}%
\]
where $R^i(w,p,\Theta)$ can be obtain from (3.22)(3.23). Indeed,
direction calculations from (3.22) show
\[
\begin{array}{ll}
-R^\alpha\eta=-\nu[-2\varepsilon^{-1}\frac{\partial^2w^\alpha}{%
\partial\xi\partial x^\beta}+2(r\varepsilon)^{-2}r^2\Theta_\beta\frac{%
\partial^2w^\alpha}{\partial(\xi)^2} -2r\varepsilon\delta_{2\alpha}\frac{%
\partial w^3}{\partial x^\beta}]\eta_\beta &  \\
\quad-\nu[-(r\varepsilon)^{-1}(\delta_{\alpha\lambda}\eta_2+2\delta_{2%
\alpha}\eta_\lambda +r\delta_{\alpha\lambda}\widetilde{\Delta}\eta)\frac{%
\partial w^\lambda}{\partial\xi} -2\varepsilon\delta_{2\alpha}a\eta_2w^3-4%
\varepsilon\delta_{2\alpha}r^2\Theta_\beta\Theta_2w^3\eta_\beta &  \\
\quad+ \delta_{2\alpha}(-2\eta_2\Theta_\sigma-2\Theta_2\eta_\sigma
-2r^2(\Theta_2\eta_\sigma+\eta_2\Theta_\sigma)|\widetilde{\nabla}\Theta|^2
-2a_{2\sigma}\Theta_\beta\eta_\beta
-r\eta_\lambda\Theta_{\lambda\sigma}-r\Theta_\lambda\eta_{\lambda\sigma})w^%
\sigma)] &  \\
\quad-[r\delta_{2\alpha}((\Pi(w,\Theta)
+2\omega)w^\beta+w^\beta\Pi(w,\Theta)) +\varepsilon^{-1}\delta_{\alpha\beta}%
\frac{\partial p}{\partial\xi}]\eta_\beta &  \\
\quad=\nu
(-r\delta_{2\alpha}w^\sigma\Theta_\lambda\eta_{\lambda\sigma}+%
\varepsilon^{-1}\frac{\partial w^\alpha}{\partial\xi}\widetilde{\Delta}\eta)+%
\nu[2\varepsilon^{-1}\frac{\partial^2w^\alpha}{\partial\xi\partial x^\beta}%
-2(r\varepsilon)^{-2}r^2\Theta_\beta\frac{\partial^2w^\alpha}{\partial(\xi)^2%
} +2r\varepsilon\delta_{2\alpha}\frac{\partial w^3}{\partial x^\beta}]%
\eta_\beta &  \\
\quad+\nu[(r\varepsilon)^{-1}(\delta_{\alpha\lambda}\delta_{2\beta}+2%
\delta_{2\alpha}\delta_{\lambda\beta}) \frac{\partial
w^\lambda}{\partial\xi}
+2\varepsilon\delta_{2\alpha}w^3(a\delta_{2\beta}+2r^2\Theta_\beta\Theta_2)]
\eta_\beta &  \\
\quad+
\nu\delta_{2\alpha}(2\delta_{2\beta}\Theta_\sigma+2\Theta_2\delta_{\beta%
\sigma} +2r^2(\Theta_2\delta_{\sigma\beta}+\delta_{2\beta}\Theta_\sigma)%
|\widetilde{\nabla}\Theta|^2 -2a_{2\sigma}\Theta_\beta
-r\Theta_{\beta\sigma})w^\sigma)\eta_\beta &  \\
\quad-[2r\delta_{2\alpha}(\Pi(w,\Theta) +\omega)w^\beta
+\varepsilon^{-1}\delta_{\alpha\beta}\frac{\partial p}{\partial\xi}%
]\eta_\beta &  \\
&
\end{array}%
\]

Taking into accounts of the homogenous boundary conditions for the
$\eta$ and integrating by parts and using Green formula and
incompressible condition
\[
\nu\varepsilon^{-1}\widetilde{\Delta}\frac{\partial w^\alpha}{\partial\xi}%
-2\nu\varepsilon^{-1}\frac{\partial^3w^\alpha}{\partial
x^\beta\partial
x^\beta\partial\xi}=-\nu\varepsilon^{-1}\widetilde{\Delta}\frac{\partial
w^\alpha}{\partial\xi},
\]
we claim
\[
\begin{array}{ll}
R^\alpha(w,p,\Theta) &
=\nu\delta_{2\alpha}\frac{\partial^2}{\partial
x^\lambda\partial x^\sigma} (rw^\sigma\Theta_\lambda)+\nu\varepsilon^{-1}%
\widetilde{\Delta}\frac{\partial w^\alpha}{\partial\xi}+\nu\frac{\partial}{%
\partial{x^\beta}} [2(\varepsilon)^{-1}\frac{\partial^2w^\alpha%
}{\partial\xi\partial x^\beta}-2\varepsilon^{-2}\Theta_\beta\frac{\partial^2w^\alpha%
}{\partial\xi^2} +2r\varepsilon\delta_{2\alpha}\frac{\partial
w^3}{\partial
x^\beta}] \\
& +\nu\frac{\partial}{\partial{x^\beta}} [(r\varepsilon)^{-1}(\delta_{\alpha%
\lambda}\delta_{2\beta}+2\delta_{2\alpha}\delta_{\lambda\beta}) \frac{%
\partial w^\lambda}{\partial\xi} +2\varepsilon\delta_{2\alpha}w^3(a\delta_{2%
\beta}+2r^2\Theta_\beta\Theta_2)] \\
& + \nu\delta_{2\alpha}\frac{\partial}{\partial x^\beta}[(2\delta_{2\beta}%
\Theta_\sigma+2\Theta_2\delta_{\beta\sigma}
+2r^2(\Theta_2\delta_{\sigma\beta}+\delta_{2\beta}\Theta_\sigma)|\widetilde{%
\nabla}\Theta|^2 -2a_{2\sigma}\Theta_\beta
-r\Theta_{\beta\sigma})w^\sigma]
\\
& -\frac{\partial}{\partial
x^\beta}[2r\delta_{2\alpha}(\Pi(w,\Theta)
+\omega)w^\beta +\varepsilon^{-1}\delta_{\alpha\beta}\frac{\partial p}{%
\partial\xi}] \\
&
\end{array}%
\]
This is (3.25). Next we consider (3.26). Indeed,
\[
\begin{array}{ll}
-R^3\eta & =-\nu[-2\varepsilon^{-1}\eta_\beta\frac{\partial^2 w^3}{%
\partial\xi\partial x^\beta}+(r\varepsilon)^{-2}2r^2\Theta_\beta\eta_\beta%
\frac{\partial^2 w^3}{\partial\xi^2} +\varepsilon^{-2}r^{-3}2r^2\Theta_\beta%
\eta_\beta \frac{\partial w^2}{\partial\xi} \\
& +2(r\varepsilon)^{-1}(\delta_{2\beta}\eta_\lambda +r\eta_{\beta\lambda})%
\frac{\partial w^\lambda}{\partial x^\beta} \\
& -(r\varepsilon)^{-1}(\eta_2+r\widetilde{\Delta}\eta)\frac{\partial w^3}{%
\partial \xi}+(r\varepsilon)^{-1}(r^{-1}\delta_{2\alpha}\eta_2+3\eta_{2%
\alpha} +r\partial_\alpha\widetilde{\Delta}\eta)w^\alpha]
+\varepsilon^{-1}w^\beta w^\lambda\eta_{\beta\lambda} \\
& +(r\varepsilon)^{-1}w^\beta\eta_\beta(2w^2
+r^2\Theta_2\Pi(w,\Theta)-2\omega(x^2)^2\Theta_2) \\
& +(r\varepsilon)^{-1}\Pi(w,\Theta)(
r^2\eta_2\Pi(w,\Theta)+r^2\Theta_2w^\beta\eta_\beta-2\omega(x^2)^2\eta_2)
-\varepsilon^{-1}\eta_\beta\partial_\beta
p+(r\varepsilon)^{-1}2r^2\Theta_\beta\eta_\beta\frac{\partial
p}{\partial\xi}
\\
& =-\nu(\varepsilon)^{-1}[2\frac{\partial w^\lambda}{\partial x^\beta}%
\eta_{\beta\lambda}-\frac{\partial w^3}{\partial\xi}\widetilde{\Delta}%
\eta+3r^{-1}w^\alpha\eta_{2\alpha} +w^\alpha\partial_\alpha\widetilde{\Delta}%
\eta ]+\varepsilon^{-1}w^\beta w^\lambda\eta_{\beta\lambda} \\
& +\nu\varepsilon^{-1}[2\frac{\partial^2 w^3}{\partial x^\beta\partial\xi}%
-2\varepsilon^{-1}\Theta_\beta\frac{\partial^2 w^3}{\partial\xi^2}%
-2(r\varepsilon)^{-1}\Theta_\beta\frac{\partial w^2}{\partial\xi}-2r^{-1}%
\frac{\partial w^\beta}{\partial x^2}+r^{-1}\delta_{2\beta}\frac{\partial w^3%
}{\partial \xi}-r^{-2}\delta_{2\beta}w^2]\eta_\beta \\
& +(r\varepsilon)^{-1}[(2w^2+r^2\Theta_2\Pi(w,\Theta)-2\omega
r^2\Theta_2)w^\beta \\
&
+\Pi(w,v\Theta)(r^2\delta_{2\beta}\Pi(w,\Theta)+r^2\Theta_2w^\beta-2\omega
r^2\delta_{2\beta})-r\partial_\beta p+2r^2\Theta_\beta\frac{\partial p}{%
\partial\xi}]\eta_\beta,%
\end{array}%
\]
By similar manner and
\[
\partial^2_{\beta\lambda}\frac{\partial w^\lambda}{\partial x^\beta}=%
\widetilde{\Delta}\frac{\partial w^\lambda}{\partial x^\lambda}
\]
we assert
\[
\begin{array}{ll}
R^3(w,p,\Theta) &
=\nu\varepsilon^{-1}(\widetilde{\Delta}(2\frac{\partial
w^\lambda}{\partial x^\lambda}-\frac{\partial w^3}{\partial \xi}-\frac{%
\partial w^\alpha}{\partial x^\alpha})+3\partial_{2\alpha}(r^{-1}w^\alpha))+%
\varepsilon^{-1}\partial_{\beta\lambda}(w^\beta w^\lambda) \\
& +2\nu\varepsilon^{-1}\widetilde{\Delta}\frac{\partial w^3}{\partial\xi}%
-2\nu\frac{\partial}{\partial x^\beta}(\Theta_\beta\frac{\partial^2 w^3}{%
\partial \xi^2}) \\
& -\nu(\varepsilon)^{-1}\frac{\partial}{\partial x^\beta}[2(r%
\varepsilon)^{-1}\Theta_\beta\frac{\partial w^2}{\partial \xi}+2r^{-1}\frac{%
\partial w^\beta}{\partial x^2}-r^{-1}\delta_{2\beta}\frac{\partial w^3}{%
\partial \xi}+r^{-2}\delta_{2\beta}w^2] \\
& +\varepsilon^{-1}\frac{\partial}{\partial x^\beta}[rw^\beta(2\frac{w^2}{r^2%
}+2\Theta_2(\Pi(w,\Theta)-\omega))
+r\delta_{2\beta}\Pi(w,\Theta)(\Pi(w,\Theta)-2\omega] \\
& +\varepsilon^{-1}\widetilde{\Delta}p+2\varepsilon^{-1}\frac{\partial}{%
\partial x^\beta}(r\Theta_\beta\frac{\partial p}{\partial \xi}) \\
& =\nu\varepsilon^{-1}(-\widetilde{\Delta}(\frac{w^2}{r})+3\partial_{2%
\alpha}(r^{-1}w^\alpha))
+\varepsilon^{-1}\partial_{\beta\lambda}(w^\beta
w^\lambda) -2\nu\frac{\partial}{\partial x^\beta}(\Theta_\beta\frac{%
\partial^2 w^3}{\partial \xi^2}) \\
& -\nu(\varepsilon)^{-1}\frac{\partial}{\partial x^\beta}[2(r%
\varepsilon)^{-1}\Theta_\beta\frac{\partial w^2}{\partial \xi}+2r^{-1}\frac{%
\partial w^\beta}{\partial x^2}-r^{-1}\delta_{2\beta}\frac{\partial w^3}{%
\partial \xi}+r^{-2}\delta_{2\beta}w^2] \\
& +\varepsilon^{-1}\frac{\partial}{\partial x^\beta}[rw^\beta(2\frac{w^2}{r^2%
}+2\Theta_2(\Pi(w,\Theta)-\omega))
+r\delta_{2\beta}\Pi(w,\Theta)(\Pi(w,\Theta)-2\omega] \\
& +\varepsilon^{-1}\widetilde{\Delta}p+2\varepsilon^{-1}\frac{\partial}{%
\partial x^\beta}(r\Theta_\beta\frac{\partial p}{\partial \xi}) \\
&
\end{array}%
\]
This is (3.26). The proof is complete.
\end{proof}

Associational variational formulation is given by
$$
\left\{%
\begin{array}{ll}
\mbox{Find}\,(\widehat{w},\widehat{p})\,\in\,H^1_\Gamma(\Omega)\times
L^2(\Omega)\,\mbox{such that} &  \\
a(\widehat{w},v)+b(w; \widehat{w},v)+b(\widehat{w};
w,v)+2(\omega\times
\widehat{w},v)-(\widehat{p},\mbox{div} v) &  \\
\quad\quad=\int_\Omega R(w,p)vdV:=(R(w,p),v), \quad \forall\quad
v\,\in\,H^1_\Gamma(\Omega), &  \\
(q, \mbox{div}\widehat{w})=o,\quad \forall\,q\,\in\,L^2(\Omega), &
\end{array}%
\right.\eqno{(3.41)}
$$
where
$$
\begin{array}{ll}
(R(w,p),v) & =\int_\Omega
g_{ij}R^iv^jdV=\int_\Omega[a_{\alpha\beta}R^\alpha(w,p,\Theta)v^\beta \\
& +\varepsilon
r^2(R^\alpha(w,p,\Theta)v^3+R^3(w,p,\Theta)v^\alpha)+%
\varepsilon^2r^2R^3(w,p,\Theta)v^3]r\varepsilon dxd\xi,%
\end{array}%
\eqno{(3.42)}
$$
For the compressible case,
$$
\left\{%
\begin{array}{l}
\mbox{div}(\hat{w}\rho+w\widehat{\rho})=0 \\
\mbox{div}(\rho \widehat{w}^iw+\rho w^i\widehat{w}+\widehat{\rho} ww^i)+2%
\widehat{\rho}(\omega\times w)^i+2\rho (\omega\omega\times\widehat{ w}%
)^i+ag^{ij}\nabla_j(\gamma\rho^{\gamma-1}\widehat{\rho}) \\
\quad -\nabla_j(A^{ijkm}e_{km}(\hat{w}))=S^i(w,\rho) \\
\end{array}%
\right.\eqno(3.43)
$$%
where
$$
\begin{array}{ll}
S^i(w,p,\Theta)=-\partial_\beta
S^{i;\beta}(w,p;\Theta)+\partial_{\lambda\sigma}^2S^{i;(\lambda,%
\sigma)}(w,p;\Theta), &  \\
&
\end{array}%
\eqno{(3.44)}
$$%
where
$$
\begin{array}{ll}
S^{\alpha;\beta}(w,p,\Theta)=\{r\delta^\alpha_2[(\delta^\beta_\lambda
\Theta_\sigma+\delta^\beta_\sigma\Theta_\lambda)w^\lambda
w^\sigma+2\varepsilon w^3w^\beta] +2\rho
r\omega\delta_{\alpha2}w^\beta
-\varepsilon^{-1}\delta_{\alpha\beta}\frac{\partial (a\rho^\gamma)}{%
\partial\xi} &  \\
\quad+2\nu \varepsilon^{-1}[-\Theta_\alpha
g^{jm}\nabla_je_{3m}^\beta(w)-
(\varepsilon^{-2}\Theta_\alpha\Theta_\beta+\delta_{\alpha\beta}g^{33})%
\nabla_3e_{33}(w) &  \\
\quad+\varepsilon^{-1}(\Theta_\alpha\delta_{\beta\lambda}+\Theta_\lambda%
\delta_{\alpha\beta}) \nabla_\lambda
e_{33}(w)+\varepsilon^{-1}(\Theta_\alpha\delta_{\beta\lambda}+2\Theta_\beta%
\delta_{\alpha\lambda} +\Theta_\lambda\delta_{\alpha\beta})\nabla_3
e_{3\lambda}(w) &  \\
\quad-\delta_{\alpha\beta}\nabla_\lambda
e_{3\lambda}(w)-\nabla_\beta
e_{3\alpha}(w)-\nabla_3 e_{\alpha\beta}(w) &  \\
\quad-r^{-1}(\delta_{\alpha2}\delta_{\beta\lambda}+3\delta_{\alpha\lambda}%
\delta_{\beta2}) (e_{3\lambda}(w)-\varepsilon^{-1}\Theta_\lambda
e_{33}(w))]\} &  \\
S^{\alpha;(\lambda,\sigma)}(w,p,;\Theta)=-2\nu\varepsilon^{-1}[(\delta_{%
\alpha\beta}\delta_{\gamma\lambda}+\delta_{\alpha\gamma}\delta_{\beta%
\lambda}) (e_{3\gamma}(w)-\varepsilon^{-1}\Theta_\gamma e_{33}(w)) &  \\
\quad+\Theta_\alpha g^{jm}(\nabla_je^{\beta\lambda}_{3m}(w,\Theta)] &  \\
&
\end{array}%
\eqno{(3.45)}
$$
$$
\begin{array}{ll}
S^{3;\beta}(w,p,\Theta)=\{-[(r\varepsilon)^{-1}(
\delta_{2\lambda}\delta^\beta_\sigma+\delta_{2\sigma}\delta^\beta_\lambda
+r^2(\delta^\beta_2\Theta_\lambda \Theta_\sigma+
\Theta_2\delta^\beta_\sigma\Theta_\lambda+\Theta_2\Theta_\sigma\delta^\beta_%
\lambda))w^\lambda w^\sigma &  \\
\quad+ 2r(\delta^\beta_\lambda\Theta_2+\delta^\beta_2\Theta_\lambda
)w^3w^\lambda+r\varepsilon\delta^\beta_2w^3w^3] &  \\
\quad+2r\omega((w^3+\varepsilon^{-1}w^\lambda\Theta_\lambda)\delta_{\beta2}+%
\varepsilon^{-1}w^\beta\Theta_2) -\varepsilon^{-1}\nabla_\beta
p+2\varepsilon^{-2}\Theta_\beta\frac{\partial p}{\partial\xi} &  \\
\quad+2\nu g^{3k}g^{jm}\nabla_je_{km}^\beta(w)
+2\nu\varepsilon^{-1}[g^{33}(4\varepsilon^{-1}\nabla_3e_{33}(w)-\nabla_\beta
e_{33}(w))+\nabla_\lambda e_{\beta\lambda}(w) &  \\
\quad+\varepsilon^{-1}(2\delta_{\beta\lambda}\nabla_\sigma
e_{3\sigma}(w)-\nabla_\lambda e_{3\beta}(w)+\nabla_\beta
e_{3\lambda}(w))\Theta_\lambda &  \\
\quad-2\varepsilon^{-2}( 2\nabla_3e_{3\lambda}(w)+\nabla_\lambda
e_{33}(w))]-2\nu[-(r\varepsilon)^{-1}e_{2\beta}(w) &  \\
\quad+e_{3\gamma}(w)[\varepsilon^{-2}r^{-1}(\delta_{\gamma2}\Theta_\beta
-3\delta_{\beta2}\Theta_\gamma)
+r\varepsilon^{-2}(\delta_{\beta\gamma}\Theta_2-\delta_{2\beta}\Theta_%
\gamma)|\nabla\Theta|^2] &  \\
\quad+e_{33}(w)(r\varepsilon)^{-3}\delta_{2\beta}(4+3r^2|\nabla\Theta|^2)]\}
&  \\
S^{3;(\lambda,\sigma)}(w,p,;\Theta)=\{2\nu[\varepsilon^{-1}(\varepsilon^{-2}%
\Theta_\lambda\Theta_\sigma +\delta_{\lambda\sigma}g^{33})e_{33}(w)
-\varepsilon^{-2}(\Theta_\lambda\delta_{\gamma\sigma}
+\Theta_\gamma\delta_{\lambda\sigma})e_{3\gamma}(w)] &  \\
\quad-\varepsilon^{-2}w^\lambda w^\sigma +2\nu
g^{3k}g^{jm}\nabla_je^{\lambda\sigma}_{km}(w,\Theta)\}.
\end{array}%
\eqno{(3.46)}
$$
%%%%%%%%%%%%%%%%%%%%%%%%%%%%%
\section{A Principle of Geometric Design for Blade's Surface}\label{sec-4}

We consider naturedly to choice  global dissipative energy to be
object functional for geometric design of the blade  surface. The
global dissipative energy functional is given by
$$
\left\{%
\begin{array}{ll}
\Phi(w,v) & =A^{ijkl}e_{kl}(w)e_{ij}(v), \\
J(S) & =1/2\int\int\int\limits_{\Omega_\varepsilon}
\Phi(w(S),w(S))\mbox{d}
V =\int\int\int\limits_{\Omega_\varepsilon}A^{ijkl}e_{kl}(w)e_{ij}(w)\sqrt{g}%
dxd\xi, \\
A^{ijkl} & =\lambda g^{ij}g^{kl}+\mu(g^{ik}g^{jl}+g^{il}g^{jk}),%
\end{array}%
\right.\eqno{(4.1)}
$$
where  $\Omega=D\times[-1,1]$ and $\Omega_\varepsilon$ is flow
passage in
Turbomachinery bounded by  $\Gamma_t\cup\Gamma_b\cup\Gamma_{in}\cup%
\Gamma_{out}\cup  S_+\cup S_-$. We proposal a principle for
geometric design of the blade:
$$
\left\{%
\begin{array}{ll}
\mbox{Find a surface}\, \Im\,\, \mbox{of the blade such that} &  \\
J(\Im)=\inf\limits_{S\in \mathcal{F}}J(S), &
\end{array}%
\right.\eqno(4.2)
$$
where $\mathcal{F}$ denotes a set of regular surfaces spanning on a
given Jordannian curve $C\in\,R^3$. The $\Im$ is called a "general
minimal surface which achieves a minimum of object dissipative
energy functional. In other words, from mathematical point of view,
this minimum problem of geometric sharp of the surface of the blade
is a general minimal surface problem.

Note that (4.2) is also an optimal control problem with distribute
parameters, control variable is the surface of the blade and the
Navier-Stokes equations are the state equations of this control
problem.

In sequel, we prove that equation (2.7) is the Eular-Lagrange
equations of optimal problem and Navier-Stokes equations are its
state equations.

In order to investigate optimal control problem (4.2) we should
consider the object functional $J$ in a fixed domain in new
coordinate system.

Let
$\widehat{w}\eta=\frac{\mathcal{D}w}{\mathcal{D}\Theta}(\Theta)\eta$
denote Gateaux  derivative of $w$ at $\Theta$ in direction $\eta\in
V(\Omega)$.Then according to (3.20) , we  have following lemma

\begin{lem}\label{lem4}
 Assume that $(w,p)$ is a solution of
Navier-Stokes equations (3.1)(3.6) associated with $\Theta\in\,H^1(D)$ which define a mapping  $%
(w(\Theta),p(\Theta)$ from $W(D)=H^1_0(D)\cap H^2(D)$ to $%
H^1(\Omega)^3\times L^2(\Omega)$  :$\Theta\Rightarrow (w(\Theta),p(\Theta))$%
. Then $(w(\Theta),p(\Theta))$ is Gateaux differential in  $H^1(D)$
with
respect to all direction  $\eta\in\,W(D)$. The Gateaux derivatives $\widehat{w%
}\eta=\frac{\mathcal{D}w}{\mathcal{D}\Theta}\eta,
\widehat{p}\eta=\frac{ \mathcal{D}p}{\mathcal{D}\Theta}\eta $ are a
solution of (3.37).  Then the strain rate tensor  $e_{ij}(w)$
defined by (3.2) possesses Gateaux derivatives
$\frac{\mathcal{D}e_{ij}(w)}{\mathcal{D}\Theta}(\Theta)\eta $ at any
point $\Theta\in\,H^2(D)$ along every direction $\eta\in\,W(D)$, and
$$
\frac{\mathcal{D}e_{ij}(w)}{\mathcal{D}\Theta}\eta=e_{ij}(\widehat{w})\eta
+e_{ij}^\lambda(w)\eta_\lambda+e_{ij}^{\lambda\sigma}(w)\eta_{\lambda\sigma},%
\eqno{(4.3)}
$$
where
$$
\begin{array}{ll}
e_{\alpha\beta}^\lambda(w)=\psi_{\alpha\beta}^\lambda(w)+
(\psi^{\lambda\sigma}_{\alpha\beta}(w)+\psi^{\sigma\lambda}(w))\Theta_\sigma
+\frac12r^2w^\sigma
(\delta_{\alpha\lambda}\Theta_{\sigma\beta}+\delta_{\beta\lambda}\Theta_{%
\alpha\sigma}), &  \\
e^{\lambda\sigma}_{\alpha\beta}(w)=\frac12r^2w^\sigma(\Theta_\alpha\delta_{%
\beta\lambda} +\Theta_\beta\delta_{\alpha\lambda}), &  \\
e_{\alpha3}^\lambda(w)=\psi_{\alpha3}^\lambda(w)
+(\psi^{\nu\lambda}_{3\alpha}+\psi^{\lambda\nu}_{3\alpha})\Theta_\nu,%
\quad e^{\lambda\sigma}_{3\alpha}(w)=\frac12\varepsilon
r^2w^\sigma\delta_{\alpha\lambda}, &  \\
e_{33}^\lambda(w)=\psi^\lambda_{33}(w),\quad
e^{\lambda\sigma}_{33}(w)=0. &
\end{array}%
\eqno{(4.4)}
$$
\end{lem}

\begin{proof}
 From (3.20) we claim
\[
\begin{array}{ll}
\frac{\mathcal{D}e_{ij}(w)}{\mathcal{D}\Theta}\eta=e_{ij}(\widehat{w})\eta
+\varphi^\lambda_{ij}(w)\eta_\lambda+\varphi^{\lambda\sigma}_{ij}(\eta_%
\lambda\Theta_\sigma+\Theta_\lambda\eta_\sigma) +R_{ij}(w,\Theta)\eta, &  \\
R_{\alpha\beta}(w,\Theta)\eta=\frac12 r^2w^\sigma[(\delta_{\alpha\lambda}%
\Theta_{\sigma\beta}+\delta_{\beta\lambda}\Theta_{\alpha\sigma})\eta_\lambda
+(\delta_{\beta\lambda}\Theta_\alpha+\delta_{\alpha\lambda}\Theta_\beta)%
\eta_{\lambda\sigma}], &  \\
R_{3\alpha}(w,\Theta)\eta =\frac12\varepsilon
r^2w^\sigma\delta_{\alpha\lambda},\quad R_{33}(e,\Theta)=0. &
\end{array}%
\]
This leads to (4.3) and (4.4). The proof is completed.
\end{proof}

\begin{lem}\label{lem5}
Under the assumptions in Lemma \ref{lem4}, the dissipative functions
$\Phi(w)$ defined by (4.1) is Gateaux differential at $\Theta\in
H^2(D)$ along any direction $\eta\in\,W$, and
$$\frac{{\cal D}\Phi(w)}{{\cal
D}\Theta}\eta=\Phi^0(\widehat{w},w)\eta+\Phi^\lambda(w,\Theta)\eta_\lambda
+\Phi^{\lambda\sigma}(w,\Theta)\eta_{\lambda\sigma},\eqno{(4.5)}$$
where
$$\left\{\begin{array}{ll}
\Phi^0(\widehat{w},w)&=2A^{ijkl}e_{ij}(\widehat{w})e_{kl}(w),\\
\Phi^\lambda(w,\Theta)&=2A^{ijkl}e^\lambda_{ij}(w)e_{kl}(w)+M^\lambda(w,\Theta),\\
 \Phi^{\lambda\sigma}(w,\Theta)&=2A^{ijkl}e^{\lambda\sigma}_{ij}(w)e_{kl}(w),
\end{array}\right.\eqno{(4.6)}$$
$$\begin{array}{ll}
M^\lambda(w,\Theta)&=4\varepsilon^{-2}(\Theta_\lambda\delta_{\alpha\beta}+\Theta_\alpha\delta_{\beta\lambda})
e_{3\alpha}(w)e_{3\beta}(w)\\
&+4\varepsilon^{-4}r^{-2}a\Theta_\lambda e_{33}(w)e_{33}(w)
 -4\varepsilon^{-1}e_{3\alpha}(w)e_{\alpha\lambda}(w)\\
&-4\varepsilon^{-3}r^{-2}(2r^2\Theta_\alpha\Theta_\lambda+(r^2|\nabla\Theta|^2+1)\delta_{\alpha\lambda})
e_{33}(w)e_{3\alpha}(w).
\end{array}\eqno{(4.7)}$$
\end{lem}

\begin{proof}
Indeed, taking (3.24) and (4.3) into account,
$$\frac{\mathcal{D}A^{ijkl}}{\mathcal{D}\Theta}\eta
 e_{ij}(w)e_{kl}(w)=M^\lambda(w,\Theta)\eta_\lambda,\quad\eqno{(4.8)}$$
 where
$M^\lambda(w,\Theta)$ is defined by (4.7). Then
$$\begin{array}{ll}

\frac{{\cal D}\Phi(w)}{{\cal
D}\Theta}\eta&=2A^{ijkl}\frac{\mathcal{D}e_{ij}(w)}{\mathcal{D}\Theta}\eta
e_{kl}(w)
+\frac{\mathcal{D}A^{ijkl}}{\mathcal{D}\Theta}\eta e_{ij}(w)e_{kl}(w),\\
\frac{{\cal D}\Phi(w)}{{\cal
D}\Theta}\eta&=2A^{ijkl}e_{ij}(\widehat{w})e_{kl}(w)\eta+(2A^{ijkl}e^\lambda_{ij}(w)+M^\lambda(w,\Theta))\eta_\lambda\\
&+2A^{ijkl}e^{\lambda\sigma}_{ij}(w)e_{kl}(w)\eta_{\lambda\sigma},
\end{array}\eqno{(4.9)}$$ From this it yields (4.5)(4.6). The proof
is complete.
\end{proof}

\begin{thm}\label{thm5}
Under the assumptions in Lemma \ref{lem2} the object functional $J$
defined by (4.1) has a $G\hat{a}teaux$ derivative $\grad
J\equiv\dfrac{{\cal D}J}{{\cal D}\Theta}$ in every direction $\eta
\in\,W:= H^2(D)\cap H^1_0(D)$. Furthermore $\grad J$ is determined
by
   $$\begin{array}{ll}
   <\grad_\Theta(J(\Theta)),\eta>=\int\int_D[\hat{\Phi}^0(w;\hat{w})\eta
   +\hat{\Phi}^\lambda(w,\Theta)\eta_\lambda
   +\hat{\Phi}^{\lambda\sigma}(w,\Theta)\eta_{\lambda\sigma}\\
\qquad   +2\mu
r^2W^{\nu\sigma}\Theta_{\nu\lambda}\eta_{\lambda\sigma}]\varepsilon
rdx^1dx^2,
   \end{array}\eqno(4.10)$$
where
   $$\left\{\begin{array}{ll}
   \hat{\Phi}^0(w;\hat{w})&=\int^1_{-1}\Phi^0(w;\hat{w})d\xi,\qquad
\widehat{\Phi}^\lambda(w,\Theta)=\int^1_{-1}\Phi^\lambda(w,\Theta)d\xi,\\
\widehat{\Phi}^{\lambda\sigma}(w,\Theta)&=2\nu r^2(\varepsilon
W^{\sigma3}_\lambda+W^{\sigma\nu}_\lambda\Theta_\nu),\\

W^{\alpha\beta}&=\int^1_{-1}w^\alpha w^\beta d\xi,\quad
W^{\alpha\beta}_\lambda:=\int^1_{-1}w^\alpha\frac{\partial
w^\beta}{\partial x^\lambda}d\xi.
\end{array}\right.\eqno{(4.11)}$$
where $\hat{w}=\dfrac{{\cal D}w}{{\cal D}\Theta}$ are the
$G\hat{a}teaux$ of the solution $(w)$ of Navier-Stokes equations at
point $\Theta$, $\Phi^0,\Phi^\lambda$ are defined by (4.6).
\end{thm}

\begin{proof}
Indeed, taking into account of (4.5) and (4.6) we assert that
 $$\begin{array}{ll}
   <\grad_\Theta(J(\Theta)),\eta>=\int\int_D\int^1_{-1}\frac{{\cal
   D}\Phi(w,\Theta)}{{\cal D}\Theta}\eta\varepsilon rd\xi dx\\
\qquad=\int\limits_D\int^1_{-1}2\nu[\Phi^0(\widehat{w},w)\eta+\Phi^\lambda(w,\Theta)\eta_\lambda
+\Phi^{\lambda\sigma}(w,\Theta)\eta_{\lambda\sigma}]r\varepsilon
d\xi dx\\
 \qquad=\int\limits_D2\mu[   \hat{\Phi}^0(w;\hat{w})\eta
   +\hat{\Phi}^\lambda(w,\Theta)\eta_\lambda
   +\hat{\Phi}^{\lambda\sigma}(w,\Theta)\eta_{\lambda\sigma}
+2\mu
r^2W^{\nu\sigma}\Theta_{\nu\lambda}\eta_{\lambda\sigma}]\varepsilon
rdx^1dx^2,
   \end{array}\eqno(4.12)$$
where we need to show only last two terms in (4.12). In fact, from
(4.4) (4.6) and (3.24)
$$\begin{array}{ll}
\Phi^{\lambda\sigma}(w,\Theta)&=2A^{ijkl}e^{\lambda\sigma}_{ij}(w)e_{kl}(w)
=4\nu\{e^{\lambda\sigma}_{\alpha\beta}(w)e_{\alpha\beta}(w)+2(g^{33}\delta^{\alpha\beta}+\varepsilon^{-2}\Theta_\alpha
\Theta_\beta)e^{\lambda\sigma}_{3\alpha}(w)e_{3\beta}(w)\\
&+g^{33}g^{33}e^{\lambda\sigma}_{33}(w)e_{33}(w)-2\varepsilon^{-1}\Theta_\beta(e^{\lambda\sigma}_{\alpha\beta}(w)
e_{3\alpha}(w)+e_{\alpha\beta}(w)e^{\lambda\sigma}_{3\alpha}(w))\\
&-2\varepsilon^{-1}g^{33}\Theta_\alpha(e^{\lambda\sigma}_{33}(w)e_{3\alpha}(w)+e_{33}(w)e^{\lambda\sigma}_{3\alpha}(w))\\
&+\varepsilon^{-2}\Theta_\alpha\Theta_\beta(e^{\lambda\sigma}_{\alpha\beta}(w)e_{33}(w)
+e_{\alpha\beta}(w)e^{\lambda\sigma}_{33}(w))\},
\end{array}$$
In view of $e^{\lambda\sigma}_{33}(w)=0$ therefore, make
rearrangement,
$$\begin{array}{ll}
\Phi^{\lambda\sigma}(w,\Theta)&=4\nu\{e^{\lambda\sigma}_{\alpha\beta}(w)(
e_{\alpha\beta}(w)-2\varepsilon^{-1}\Theta_\beta
e_{3\alpha}(w)+\varepsilon^{-2}\Theta_\alpha\Theta_\beta
e_{33}(w))\\
&+2e^{\lambda\sigma}_{3\alpha}(w)((g^{33}\delta^{\alpha\beta}+\varepsilon^{-2}\Theta_\alpha
\Theta_\beta)e_{3\beta}(w)-\varepsilon^{-1}\Theta_\beta
e_{\alpha\beta}(w)-\varepsilon^{-1}g^{33}\Theta_\alpha
e_{33}(w))\}.\\
\end{array}$$
 Substituting (4.4) into above leads to
$$\begin{array}{ll}
\Phi^{\lambda\sigma}(w,\Theta)&=2\nu
r^2w^\sigma[(2\varepsilon^{-2}|\nabla\Theta|^2-2g^{33})\Theta_\lambda
e_{33}(w)\\
&+2\varepsilon(g^{33}\delta_{\lambda\beta}+\varepsilon^{-2}\Theta_\lambda\Theta_\beta-\varepsilon^{-2}
(\Theta_\beta\Theta_\lambda+|\nabla\Theta|^2\delta_{\beta\lambda})e_{3\beta}(w))]\\
&=4\nu \varepsilon^{-2}w^\sigma[\varepsilon
e_{3\lambda}(w)-\Theta_\lambda e_{33}(w)]\\
&=4\nu \varepsilon^{-2}w^\sigma[\frac{\varepsilon}{2}(\frac{\partial
w^\lambda}{\partial \xi}+\varepsilon^2r^2\frac{\partial
w^3}{\partial x^\lambda})+\frac12\varepsilon^2r^2(\frac{\partial
w^\nu}{\partial x^\lambda}-\frac{\partial w^3}{\partial
\xi}\delta_{\lambda\nu})\Theta_\nu\\
&-\frac12\varepsilon r^2\frac{\partial w^\nu}{\partial
\xi}\Theta_\nu\Theta_\lambda+\frac12\varepsilon^2r^2w^\nu\Theta_{\lambda\nu}],
\end{array}\eqno{(4.13)}$$

By virtue of (3.20-22) and the index $(\lambda,\sigma)$ in
$\Phi^{\lambda\sigma}\eta_{\lambda\sigma}$ being symmetry,
$$\int_{-1}^1w^\sigma \partial_\xi w^\lambda
d\xi=\int_{-1}^1\frac12\partial_\xi(w^\lambda w^\sigma)d\xi=\frac12
w^\lambda w^\sigma|_{\xi=1,-1}=0,\mbox{(by boundary conditions)}$$
we have
$$\begin{array}{ll}
\int^1_{-1}\Phi^{\lambda\sigma}(w,\Theta)d\xi&= 2\nu
r^2\int^1_{-1}w^\sigma[\varepsilon \frac{\partial w^3}{\partial
x^\lambda}+\frac{\partial w^\nu}{\partial
x^\lambda}\Theta_\nu]d\xi+2\nu
r^2W^{\nu\sigma}\Theta_{\lambda\nu}\\
&=2\nu r^2[\varepsilon
W^{\sigma3}_\lambda+W^{\sigma\nu}_\lambda\Theta_\nu+W^{\nu\sigma}\Theta_{\lambda\nu}].
\end{array}$$
 The proof is complete.
 \end{proof}

Taking integration by part of (4.10) and considering homogenous
boundary condition for $\eta\in\,W(D)$,it implies
$$<\grad_\Theta(J(\Theta)),\eta>=\int\limits_D \varepsilon
[\partial_{\lambda\sigma}(2\mu
r^3W^{\sigma\nu}\Theta_{\nu\lambda}+r\widehat{\Phi}^{\lambda\sigma}(w,\Theta))
-\partial_\lambda(r\widehat{\Phi}^\lambda(w,\Theta))
+r\widehat{\Phi}^0(w,\widehat{w})] \eta dx.$$

From above discussion we obtain directly Theorem \ref{thm1} The
Euler-Lagrange equation for the extremum $\Theta$ of $J$ is given
by:
$$\left\{\begin{array}{ll}
\frac{\partial^2}{\partial x^\lambda\partial x^\sigma}(2\nu
r^3W^{\nu\sigma}\frac{\partial^2\Theta}{\partial x^\nu\partial
x^\lambda})+
 \frac{\partial^2}{\partial x^\lambda \partial
x^\sigma}(r\hat{\Phi}^{\lambda\sigma}(w,\Theta))
-\frac{\partial}{\partial x^\lambda}(r\hat{\Phi}^\lambda(w,\Theta))
+\hat{\Phi}^0(w,\widehat{w})r=0,\\
\Theta|_\gamma=\Theta_0,\quad \frac{\partial\Theta}{\partial
n}|_\gamma=\Theta_*,
\end{array}\right.\eqno(4.14)$$
and variational formulation associated with (4.14) reads
$$\left\{\begin{array}{ll}
\mbox{Find}\Theta\in\,V(D)=\{q|q\in\,H^2(D),
q|_\gamma=\Theta_0,\frac{\partial q}{\partial n}|_\gamma=\Theta_*\}
\mbox{such that}\\
\int_D\{({\cal
K}^{\alpha\beta\lambda\sigma}\Theta_{\alpha\beta}+\hat{\Phi}^{\lambda\sigma}(w,\Theta))
\eta_{\lambda\sigma}+\widehat{\Phi}^\lambda(w,\Theta)\eta_\lambda+\widehat{\Phi}^0(w,\widehat{w})\eta\}\varepsilon
rdx,\quad \forall\,\eta\in\,H^2_0(D),
\end{array}\right.\eqno{(4.15)}$$ where
$${\cal K}^{\alpha\beta\lambda\sigma}:=2\nu
r^3W^{\alpha\sigma}\delta_{\beta\lambda}.\eqno{(4.16)}$$

\section {Second Model}\label{sec-5}

Let us choice the rate of work done by the impeller under certain
ratio of inlet and outlet pressure as minimizing functional (3.10)
$$I(\Im)=\int\int_{\Im_-\cup\Im_+}\sigma\cdot n\cdot e_\theta
r\omega dS,\eqno{(5.1)}$$ where $n$ is unite normal vector to the
surface $\Im$ and $(e_r,e_\theta ,k)$ are bases vector in rotating
cylinder coordinate system, $\sigma$ is total stress tensor for
compressible flow. Our purpose is to find a surface $\Im$ of blade
such that
$$I(\Im)=\inf_{S\in {\cal F}}I(S),\eqno{(5.2)}$$
where ${\cal F}$ denotes a set of regular surface spanning on a
given Jordannian curve $C\in R^3$. Under new coordinate system (5.1)
can be rewritten by
$$\begin{array}{ll}
I(\Im)=\int_D\{[((-p+\lambda\mbox{div} w)g^{ij}+2\mu
e_{ij}(w))n^i(e_\theta)^j\omega
r\sqrt{a}]|_{\xi=-1}\\
\qquad+[((-p+\lambda\mbox{div} w)g^{ij}+2\mu
e_{ij}(w))n^i(e_\theta)^j\omega r\sqrt{a}]|_{\xi=-1}\} r\varepsilon
dx,\end{array}\eqno{(5.3)}$$

Next we compute (5.3). To do this, we give the followings without
proof:
$$\left\{\begin{array}{ll}
n(\xi)=n^i(\xi)e_i, \quad
n^\lambda(\xi)=-r\Theta_\lambda/\sqrt{a},\quad n^3(-1)=(\varepsilon
r)^{-1} \frac{1+r^2\Theta_2^2}{\sqrt{a}},\\ n^3(1)=-(\varepsilon
r)^{-1}\frac{1+r^2\Theta_2^2}{\sqrt{a}},\quad
e_\theta=(r\varepsilon)^{-1}e_3,\quad (e_\theta)^\alpha=0,\quad
(e_\theta)^3=(r\varepsilon)^{-1}.
\end{array}\right.\eqno{(5.4)}$$

Thanks to (2.5) and (5.4) we claim that
$$\left\{\begin{array}{ll}
(-p+\lambda\mbox{div} w)g_{ij}n^i(e_\theta)^j=(-p+\lambda\mbox{div}
w)(r\varepsilon)^{-1}(g_{33}n^3+g_{3\lambda}n^\lambda)\\
\quad=(-p+\lambda\mbox{div} w)( r\varepsilon n^3+r\Theta_\alpha
n^\alpha),\quad
\mbox{for }\xi=\pm 1,\\
2\mu e_{ij}(w)n^i(e_\theta)^j=2\mu(r\varepsilon)^{-1}(e_{33}(w)n^3+e_{3\lambda}(w)n^\lambda),\\
\end{array}\right. \eqno{(5.5)}$$

Substituting (5.4)(5.5) into (5.3) and simple calculations  lead to
$$\left\{\begin{array}{ll}
I(\Im)=\int_D\{(1+r^2\Theta_2^2)(-2\widetilde{P}+2\lambda\mbox{div}
\widetilde{W})-r^2|\nabla\Theta|^2(-2P+2\lambda\mbox{div} W)\\
\quad
+2\mu(r\varepsilon)^{-1}((r\varepsilon)^{-1}(1+r^2\Theta_2^2)e_{33}(\widetilde{W})-r\Theta_\alpha
e_{3\alpha}(W))\}r\omega\varepsilon rdx,
\end{array}\right.\eqno{(5.6)}$$
where we use notation
$$\left\{\begin{array}{ll}
w|_{\xi=-1}=w(-1),\quad w|_{\xi=1}=w(1),\quad W=(w(1)+w(-1))/2,\quad
\widetilde{W}=(w(-1)-w(1))/2,\\
P=(p(1)+p(-1))/2,\quad \widetilde{P}=(p(-1)-p(1))/2,\\
W_\Theta:=\frac{{\cal D}W}{{\cal D}\Theta},\quad
P_\Theta:=\frac{{\cal D}P}{{\cal D}\Theta},\quad
\widetilde{W}_\Theta:=\frac{{\cal D}\widetilde{W}}{{\cal
D}\Theta},\quad \widetilde{P}_\Theta:=\frac{{\cal
D}\widetilde{P}}{{\cal D}\Theta},
\end{array}\right.\eqno{(5.7)}$$
Hence
$$\widetilde{P}-P=-p(-1), \widetilde{W}-W=-w(-1),\eqno{(5.8)}$$
and in view of Lemma \ref{lem2}, the gradient of $I$ along any
direction $\eta\in H^1_0(D)$ is given by

\begin{thm}\label{thm6} Under the assumptions in Lemma \ref{lem2} the object functional $I$
defined by (5.1) has a $G\hat{a}teaux$ derivative $\grad
I\equiv\dfrac{{\cal D}I}{{\cal D}\Theta}$ in every direction $\eta
\in\, H^1_0(D)$. Furthermore $\grad I$ is determined by
$$\begin{array}{ll}
<\grad
I(\Im),\eta>=\int_D\{[\Psi_0(w,p,\Theta)\eta+\Psi^\lambda(w,p,\Theta)\eta_\lambda
+\Psi^{\lambda\sigma}(w,p,\Theta)\eta_{\lambda\sigma}]r\omega\varepsilon
r\}dx,\end{array}\eqno{(5.9)}$$where
$$\left\{\begin{array}{ll}
\Psi_0(w,p,\Theta):=
(1+r^2\Theta_2^2)(-2\widetilde{P}_\Theta+2\lambda\mbox{div}
\widetilde{W}_\Theta)-r^2|\nabla\Theta|^2(-2P_\Theta+2\lambda\mbox{div} W_\Theta)\\
\quad
+2\mu(r\varepsilon)^{-1}((r\varepsilon)^{-1}(1+r^2\Theta_2^2)e_{33}(\widetilde{W_\Theta})-r\Theta_\alpha
e_{3\alpha}(W_\Theta)),\\
\Psi^\lambda(w,p,\Theta)=\Psi_0^\lambda(w,p)
+\Psi_\nu^\lambda(w,p)\Theta_\nu+
\Psi_{\nu\mu}^\lambda(w,p)\Theta_\nu\Theta_\mu+\Psi_*^\lambda(w,p,\Theta),\\

\Psi^{\lambda\sigma}(w,p,\Theta):=-2\mu
\varepsilon^{-1}\Theta_\alpha e_{3\alpha}^{\lambda\sigma}(W)=- \mu
r^2W^\sigma\Theta_\lambda,
\end{array}\right.\eqno{(5.10)}$$

$$\left\{\begin{array}{ll}
\Psi_0^\lambda(w,p)&=2\mu((r\varepsilon)^{-2}\psi^\lambda_{33}(\widetilde{W})-
\varepsilon^{-1}\varphi_{3\lambda}(W)),\\
\Psi_\nu^\lambda(w,p)&=[-2r^2(-2P+2\lambda\mbox{div}
W)\delta_{\lambda\nu} +(2r^2(-\widetilde{P}+2\lambda\mbox{div}
\widetilde{W})\\
&+2\varepsilon^{-2}\varphi_{33}(\widetilde{W}))\delta_{2\lambda}\delta_{2\nu}
-2\mu\varepsilon^{-1}(\psi^\lambda_{3\nu}(W)+\psi^\nu_{3\lambda}(W))],\\
\Psi_{\nu\mu}^\lambda(w,p)&=2\mu
\varepsilon^{-2}[\psi^\lambda_{33}(\widetilde{W})\delta_{2\nu}+2\psi^\nu_{33}(\widetilde{W})
\delta_{2\lambda}]\delta_{2\mu}\\
&-2\mu\varepsilon^{-1}[\psi^{\nu\lambda}_{3\mu}(W)
+\psi^{\lambda\nu}_{3\mu}(W)+\psi^{\nu\mu}_{3\lambda}(W)],\\
\Psi_*^\lambda(w,p,\Theta)&=-\mu r^2W^\sigma\Theta_{\lambda\sigma},\\
\end{array}\right.\eqno{(5.11)}$$
\end{thm}

\begin{proof}
Firstly,   by virtue of (3.21) we assert
$$\frac{{\cal D}\mbox{div} W}{{\cal D}\Theta}=\mbox{div} \frac{{\cal D}W}{{\cal
D}\Theta}.$$ In addition, thanks to (4.7) and denote
$$W_\Theta:=\frac{{\cal D}W}{{\cal D}\Theta},\quad \widetilde{W}_\Theta:
=\frac{{\cal D}\widetilde{W}}{{\cal D}\Theta},\quad
P_\Theta:=\frac{{\cal D}P}{{\cal D}\Theta},\quad
\widetilde{P}_\Theta:=\frac{{\cal D}\widetilde{P}}{{\cal
D}\Theta},$$ (5.6) shows
$$\begin{array}{ll}
<\grad_\Theta
I,\eta>=\int_D\{2r^2\Theta_2\eta_2(-2\widetilde{P}+2\lambda\mbox{div}
\widetilde{W})+(1+r^2\Theta_2)(-2\widetilde{P}_\Theta+2\lambda\mbox{div}
\widetilde{W}_\Theta)\eta\\
\quad-2r^2\Theta_\lambda\eta_\lambda(-2P+2\lambda\mbox{div}
W)-r^2|\widetilde{\nabla}\Theta|^2(-2P_\Theta+2\lambda\mbox{div}
W_\Theta)\eta\\
\quad+2\mu(r\varepsilon)^{-1}[(r\varepsilon)^{-1}2r^2\Theta_2\eta_2
e_{33}(\widetilde{W})+(r\varepsilon)^{-1}(1+r^2\Theta_2^2)\frac{{\cal
D}e_{33}}{{\cal D}\Theta}(\widetilde{W})\eta-r\eta_\alpha
e_{3\alpha}(W)-r\Theta_\alpha\frac{{\cal D}e_{3\alpha}}{{\cal
D}\Theta}(\widetilde{W})\eta ]
\end{array}$$ In view of (4.7)(4.8) and rewriting integrated we
claim the first of (5.10). Moreover
$$\begin{array}{ll}
\Psi^\lambda(w,p,\Theta)&=2r^2\Theta_2\delta_{2\lambda}(-\widetilde{P}+2\lambda\mbox{div}
\widetilde{W})-2r^2\Theta_\lambda(-2P+2\lambda\mbox{div} W)\\
&+2\mu(r\varepsilon)^{-1}[(r\varepsilon)^{-1}(1+r^2\Theta_2)e^\lambda_{33}(\widetilde{W})
+2r^2(r\varepsilon)^{-1}\delta_{2\lambda}\Theta_2e_{33}(\widetilde{W})\\
&-r\Theta_\alpha
e^\lambda_{3\alpha}(W)-re_{3\lambda}(W)],\\
\end{array}\eqno{(5.12)}$$
$$\Psi^{\lambda\sigma}(w,\Theta)=-2\mu\varepsilon^{-1}\Theta_\alpha
e^{\lambda\sigma}_{3\alpha}(W)=- \mu r^2W^\sigma\Theta_\lambda.$$
 Thanks to
$$e_{3\lambda}(w)=\psi_{3\lambda}(w)+\psi^\nu_{3\lambda}(w)\Theta_\nu
+\psi^{\nu\mu}\Theta_\nu\Theta_\mu+\frac12\varepsilon
r^2w^\sigma\Theta_{\lambda\sigma}$$ and Taking (4.8) and
(3.16)(3.17) into account, simple calculation from (5.12) shows
$$\begin{array}{ll}
\Psi^\lambda(w,p,\Theta)&=2\mu(r\varepsilon)^{-2}\psi^\lambda_{33}(W)-2\mu\varepsilon^{-1}\varphi_{3\lambda}(W)\\
&+[(2r^2(-\widetilde{P}+2\lambda\mbox{div}
\widetilde{W})+4\mu\varepsilon^{-2}\varphi_{33}(\widetilde{W}))\delta_{2\lambda}\delta_{2\nu}\\
&-2r^2(-2P+2\lambda\mbox{div}
W)\delta_{\lambda\nu}-2\mu\varepsilon^{-1}(\psi^\lambda_{3\nu}(W)+\psi^\nu_{3\lambda}(W))]\Theta_\nu\\
&+[2\mu\varepsilon^{-2}(\psi^\lambda_{33}(W)\delta_{2\nu}+2\psi^\nu_{33}(W)\delta_{2\lambda})\delta_{2\mu}\\
&-2\mu\varepsilon^{-1}(\psi^{\nu\lambda}_{3\mu}(W)+\psi^{\lambda\nu}_{3\mu}(W)
+\psi^{\nu\mu}_{3\lambda}(W))]\Theta_\nu\Theta_\mu\\
&-2\mu(\varepsilon r)^{-1}r\frac12\varepsilon r
W^\sigma\Theta_{\lambda\sigma},
\end{array}\eqno{(5.13)}$$
This leads to second of (5.10).

 Integrating by parts in (5.9)  yields
$$\begin{array}{ll}
<Grad_\Theta I, \eta>=\int_D[\varepsilon
r^2\omega\Psi_0(w,p,\Theta)\\
\quad-\frac{\partial}{\partial x^\lambda}(\varepsilon r^2\omega
\Psi^\lambda(w,p,\Theta))+\frac{\partial^2}{\partial x^\lambda
\partial x^\sigma}(\varepsilon\omega
r^2\Psi^{\lambda\sigma}(w,p,\Theta))]\eta dx
\end{array}\eqno{(5.14)}$$
Hence the stationary point of minimum problem should satisfies
following Eular-Lagrange equation
$$\begin{array}{ll}
\varepsilon r^2\omega\Psi_0(w,p,\Theta) -\frac{\partial}{\partial
x^\lambda}(\varepsilon r^2\omega
\Psi^\lambda(w,p,\Theta))+\frac{\partial^2}{\partial x^\lambda
\partial x^\sigma}(\varepsilon\omega
r^2\Psi^{\lambda\sigma}(w,p,\Theta))=0
\end{array}\eqno{(5.15)}$$
However,
$$\begin{array}{ll}
\frac{\partial^2}{\partial x^\lambda
\partial x^\sigma}(\varepsilon\omega
r^2\Psi^{\lambda\sigma}(w,p,\Theta))=-\frac{\partial^2}{\partial
x^\lambda
\partial x^\sigma}(\mu\omega\varepsilon
r^4W^\sigma\Theta_\lambda)\\
\quad=-\frac{\partial^2}{\partial x^\lambda
\partial x^\sigma}(\mu\omega\varepsilon
r^4W^\sigma)\Theta_\lambda-\mu\varepsilon\omega r^4
W^\sigma\frac{\partial^2}{\partial x^\lambda
\partial x^\sigma}(\Theta_\lambda)-\frac{\partial}{\partial x^\lambda
}(\mu\omega\varepsilon r^4W^\sigma)\frac{\partial}{
\partial x^\sigma}(\Theta_\lambda)\\
\quad-\frac{\partial}{\partial x^\sigma }(\mu\omega\varepsilon
r^4W^\sigma)\frac{\partial}{
\partial x^\lambda}(\Theta_\lambda)
\end{array}$$
In addition
$$\begin{array}{ll}
-\frac{\partial}{\partial x^\lambda
}(\Psi^\lambda_*(w,p,\Theta))=\frac{\partial}{\partial x^\lambda
}(\varepsilon\omega\mu r^4W^\sigma\Theta_\lambda)
=\frac{\partial}{\partial x^\lambda}(\varepsilon\omega\mu
r^4W^\sigma)\Theta_{\sigma\lambda}+\mu\varepsilon\omega r^4
W^\sigma\Theta_{\sigma\lambda\lambda},
\end{array}$$ Therefore
$$\begin{array}{ll}
-\frac{\partial}{\partial x^\lambda
}(\Psi^\lambda_*(w,p,\Theta))+\frac{\partial^2}{\partial x^\lambda
\partial x^\sigma}(\varepsilon\omega
r^2\Psi^{\lambda\sigma}(w,p,\Theta)) = -\frac{\partial^2}{\partial
x^\lambda
\partial x^\sigma}(\mu\omega\varepsilon
r^4W^\sigma)\Theta_\lambda-\frac{\partial}{\partial x^\sigma
}(\mu\omega\varepsilon r^4W^\sigma)\widetilde{\Delta}\Theta,\\
\end{array}$$where
$$\widetilde{\Delta}\Theta=\Theta_{11}+\Theta_{22}.
$$ On the other hand,
$$\begin{array}{ll}
-\partial_\lambda(\Psi^\lambda_0+\Psi^\lambda_\nu\Theta_\nu+\Psi^\lambda_{\mu\nu}\Theta_\nu\Theta_\mu)
=-(\partial_\lambda\Psi^\lambda_0+\partial_\lambda\Psi^\lambda_\nu\Theta_\nu+\partial_\lambda
\Psi^\lambda_{\nu\mu}\Theta_\nu\Theta_\mu)
+(\Psi^\lambda_\nu+(\Psi^\lambda_{\nu\mu}+\Psi^\lambda_{\mu\nu})\Theta_\mu)\Theta_{\lambda\nu}
\end{array}$$
To sum up and by simply calculation, we obtain Theorem \ref{thm2}:
$$\left\{\begin{array}{ll}
-(K_0(w)\widetilde{\Delta}\Theta+K^{\lambda\nu}(w,\Theta)\Theta_{\nu\lambda})
+F^{\nu\mu}(w)\Theta_\nu\Theta_\mu
+F^\lambda(w)\Theta_\lambda+F_0(w,\Theta)=0,\\
\Theta|_\gamma=\Theta_0,\quad \gamma=\partial D,
\end{array}\right.\eqno{(5.16)}$$
$$\left\{\begin{array}{ll}
K_0(w)=\mu\varepsilon\omega \frac{\partial(r^4W^\sigma)}{\partial x^\sigma},\\
K^{\lambda\nu}(w)=\omega\varepsilon
r^2(\Psi^\lambda_\nu(w)+(\Psi^\lambda_{\nu\mu}(w)+\Psi^\lambda_{\mu\nu}(w))\Theta_\mu),\\
F^{\nu\mu}(w)=-\omega\varepsilon r(r\partial_\lambda
\Psi^\lambda_{\nu\mu}(w)+2\delta_{2\lambda}\Psi^\lambda_{\nu\mu}(w)),\\
F^\lambda(w)=-(\mu\varepsilon\omega\partial^2_{\nu\sigma}(
r^4W^\sigma)+\omega\varepsilon
r(r\partial_\lambda\Psi^\lambda_\nu(w)+2\delta_{2\lambda}\Psi^\lambda_\nu(w)),\\
F_0(w,\Theta)=\omega\varepsilon
r^2\Psi_0(w,\Theta)-\omega\varepsilon
r(r\partial_\lambda\Psi^\lambda_0(w)+2\delta_{2\lambda}\Psi^\lambda_0(w)),
\end{array}\right.\eqno{(5.17)}$$ where $\Psi_0,\Psi^\lambda_\nu $ and $\Psi^\lambda_{\nu\mu}$
are defined by (5.10) and (5.11) respectively.  Variational
formulation associated with (5.16) is given by
$$\left\{\begin{array}{ll}
\mbox{Find}\,\Theta\,\in\,V=\{q|q\in\,H^1(D),q|_\gamma=\Theta_0\}\quad\mbox{such
that}\quad\\
\int_D\{[\Psi_0(w,p,\Theta)\eta+\Psi^\lambda(w,p,\Theta)\eta_\lambda
+\Psi^{\lambda\sigma}(w,p,\Theta)\eta_{\lambda\sigma}]r\omega\varepsilon
r\}dx,
\end{array}\right.\eqno{(5.18)}$$
where $\Psi_0(w,p,\Theta),\Psi^\lambda(w,p,\Theta),
\Psi^{\lambda\sigma}(w,p,\Theta)$ are defined by (5.10).
\end{proof}
%%%%%%%%%%%%%%%%%%%%%%%%
\section{Existence of  Solution to the Optimal Control
Problem}\label{sec-6}

In this section,  we discuss  existence of the optimal control
problem (3.1) for incompressible case. As well known the object
functional
$$J(\Theta)=\frac12\int_\Omega A^{ijkl}(\Theta)e_{ij}(w(\Theta))e_{kl}(w(\Theta))\sqrt{g}dxd\xi
=\frac12a(w(\Theta),w(\Theta)), \eqno{(6.1)}$$ is  depending upon
the existence of solution $w$ to Navier-Stokes equations. Since the
solution $w$ of Navier-Stokes equation and $A^{ijkl}$ are  the
functions of $\Theta$ , therefore $J$ is a function of $\Theta$.
However $J$ can be read as a function of $w$ :
$J(\Theta)=\widetilde{J}(w(\Theta))$. As well known that if there
exists Gateaux derivative $\frac {D J}{D\Theta}$
  of $J(\Theta)$ with respect to $\Theta$ at $\Theta^*$,
Then the minimum point $\Theta$ of (3.1) is necessary to satisfies
$$ grad_S J(\Theta)=0.$$
 and from this, if $\Theta$ is regular  then Eular-Lagrange equation for
 $\Theta$ can be obtained.

 At first,It is well known that following theorem give a sufficient
condition of the existence
\begin{thm}\label{thm7}
(Generalized Weierstrass Theorem)[13] Let X be a reflexive Banach
Space,and U a bounded and weakly closed subset of X. If the
functional J is weakly lower semi-continuous on U,then J is bounded
from below and J achieves its infirm on U.
\end{thm}
We consider functional defined in a closed convex set in a Sobolev
space
$$V(\Omega):=\{u|u\in\,H^{1,p}(\Omega),u|_{\Gamma_0}=0,\partial
\Omega=\Gamma_0 \cup\Gamma_1,meas(\Gamma_0)\neq 0\}.$$ Let
 $$\widetilde{J}(w)=\int_\Omega A^{ijkl}e_{ij}(w)e_{kl}(w)\sqrt{g}dxd\xi $$

\begin{lem}\label{lem6} $\widetilde{J}(w)$ is weak lower semi-continuous with respect
to $w$ in $H^1(\Omega)^3$.
\end{lem}
\begin{proof} We firstly need to establish the uniform positive
 definiteness of the three-dimensional tensor $A^{ijkl}$ . "Uniform" means with respect
 to all points $x\in\overline\Omega$
and to symmetric matrices of the order three $\{t_{ij}\}$. Namely,
there exist a constant $c(\Omega,\Theta,\mu)>0$ such that
 $$A^{ijkl}(x)t_{kl}t_{ij}\geq c\sum_{i,j}|t_{i,j}|^2.\eqno{(6.2)}$$
 See the proof of  Th1.8-1 in [P.G.Ciarlet 2000]. Note that
 $\Omega$ is a Lipschitz domain,$g=det(g_{ij})=\varepsilon^2r^2>0$.
Therefore, we claim
$$a(w,w)\geq \nu c(\Theta,\Theta)\sum_{i,j}\|e_{ij}(w)\|^2_{0,\Omega}$$
On the other hand,  the Korn inequality with boundary conditions in
curvilinear coordinates  (see Th1.7-4 in [P.G.Ciarlet 2000]) shows
that
 $$\sum_{i,j}||e_{ij}(v)||^2_{0,\Omega}\geq
 c(\Omega,\Theta)||v||^2_{1,\Omega},\quad \forall v\in V.\eqno{(6.3)}$$
Hence
$$ a(w,w)\geq \mu c(\Omega,\Theta)\|w\|^2_{1,\Omega}, \eqno{(6.4)}$$
We assert that $a(\cdot,\cdot)$ is a equivalent norm in
$H^1(\Omega)^3$. In other words  ,  $\widetilde{J}(\cdot)$ is a
equivalent norm in $H^1(\Omega)^3$. As well known that the norm in a
Hilbert space , as a functional, is weakly lower semi-continuous. we
assert $\widetilde{J}$ is weakly lower semi-continuous with respect
to $w$. Proof is complete.
\end{proof}

 By virtue of Lemma \ref{lem3} we directly to obtain
\begin{lem}\label{lem7}
If  the function $w(\Theta)$ of Navier-Stokes
 equations
satisfies the following:
$${\bf \mbox{Assumption}\quad P}:\Theta_n\rightharpoonup \Theta_0(weakly)\Rightarrow
w_n=w(\Theta_n)\rightharpoonup w_0=w(\Theta_0)(weakly).$$ Then
functional $J(\Theta)$ is weak lower semi-continuous with respect to
$\Theta$. \end{lem}

Finally, we have
\begin{thm}\label{thm8}Assume that the solutions
$(w,p)$ of Navier-Stokes
 equations with mixed
 boundary condition are weakly continuous with respect to $\Theta$.Then
 there exists a two dimensional surface $S$ defined by a smooth mapping
  $$\Theta: D \longrightarrow {\cal W}\equiv H^2(D)\cap (\Theta\in H^1(D),
  \Theta|_{\partial D}=\Theta^*,)$$
  such that $J(\Theta)$   achieves its infirmum at $\{\Theta
 ,w(\Theta)\}$.
\end{thm}

 Next we consider the existence of the solutions for
Navier-Stokes equations. In deed flow's domain is unbounded domain.
In section 3 we make artificial boundary, inflow boundary
$\Gamma_{in}$ and outflow boundary $\Gamma_{out}$, and impose
natural boundary conditions (3.6). We also can impose the pressures
$$p|_{\Gamma_{in}}=p_{in},\quad p|_{\Gamma_{out}}=p_{out},$$
or flux
$$\int_{\Gamma_{in}}\rho w\cdot nd\Gamma=Q,\quad \int_{\Gamma_{out}}\rho w\cdot
nd\Gamma=Q,$$

Let us consider energy inequality. Owing to
$$ (2\omega\times w,w)=0,$$ Hence
moment equations (3.11) show
$$a(w,w)+b(w,w,w)=(f,w),$$
However,
$$\begin{array}{ll}
b(w,w,w)=\int_\Omega w^j\nabla_jw^ig_{ik}w^k\sqrt{g}dxd\xi
=\int_\Omega(\nabla_j(w^jw^i)-w^i\div
w)g_{ik}w^k\sqrt{g}dxd\xi\\
\quad=\int_\Omega(\div(|w|^2w)-g_{ik}w^iw^j\nabla_jw^k)\sqrt{g}dxd\xi
= \int_{\Gamma_1}|w|^2w\cdot nd\Gamma-b(w,w,w) \end{array}$$
$$b(w,w,w)=\frac12\int_{\Gamma_1}|w|^2w\cdot nd\Gamma,\eqno{(6.5)}$$Here we
denote
$$|w|^2=g_{ik}w^iw^k,\quad \Gamma_1=\Gamma_{in}\cup\Gamma_{out}.\eqno{(6.6)}$$
The flux of inflow and outflow kinetic energy are respectively give
by
$$ K_{in}(w)=\int_{\Gamma_{in}}|w|^2w\cdot nd\Gamma,\quad
K_{out}(w)=\int_{\Gamma_{out}}|w|^2w\cdot nd\Gamma$$ where $w\cdot
n=g_{ij}w^in^j$ and $n$ is outward normal unite vector to inflow or
out flow boundary. Therefore (6.1) shows that
$$b(w,w,w)=K_{out}(w)-K_{in}(w),\eqno{(6.7)}$$
It is obvious that we can not make say $K_{out}(w)<\infty$ .

\begin{thm}\label{thm9}
 Suppose that the exterior force $f$ and normal stress $g$ of
inflow and out flow boundaries
$\Gamma_1=\Gamma_{in}\cup\Gamma_{out}$ satisfy
$$\|F\|_*:=\|f\|_{0,\Omega}+\|g_{in}\|_{-1/2,\Gamma_{in}}+\|g_{out}\|_{-1/2,\Gamma_{out}}\leq
\nu^2/(4c_1c_2),\eqno{(6.8)}$$ and the mapping $\Theta$ defined in
(2.1) is $C^2(D)$-function satisfied by
$$\Theta\in\,
C^2(D,K_0)\equiv\{\theta\in\,C^2(D),\|\theta_\lambda\|\leq K_0;
,\|\theta_{\lambda\sigma}\| \leq K_0\},\eqno{(6.9)}$$ where $K_0$ is
a constant.
 Then there exists a smooth solution of the
variational problem (3.11) satisfying
$$\|\nabla w\|_{0,\Omega}\leq
\frac{c(\Omega,\Theta)\mu}{2c_1}[1-\sqrt{(1-\frac{4c_1c_2\|F\|_*}{(c(\Omega,\Theta)\mu)^2})}],\eqno{(6.10)}$$
where $c_1,c_2$ are constants depending of $(\Omega, K_0)$ defined
by (6.14)(6.15).
 \end{thm}

\begin{proof}
To prove the theorem for a steady Navier-Stokes problem, it is
convenient to construct the solution as a limit of Galerkin
approximations in terms of the eigenfunctions of the corresponding
stead Stokes problem. This use of the Stokes eigenfunctions (see in
Heywood and Rannacher and Turek[20] and Glowinski[12]). Galerkin
equations are a system of algebraic equations for constant unknowns
and the Galerkin approximation solution $w$ is a solution of the
finite dimensional problem
$$ a(w,v)+2(\omega\times w,v)+b(w,w,v)=<F,v>, \forall
v\in\,V_m:=\mbox{span}\{\varphi_1,\varphi_2\cdots,\varphi_m\}\eqno{(6.11)}$$
where $\varphi_i, i=1,2,\cdots m$ are eigenfunctions of
corresponding Stokes operator. Let $S_\rho$ denote the spheres in
$V_m$ satisfying inequality (6.10). Assume that $w_*\in\,S_\rho$, we
seeking $w$ such that
$$a(w,v)+2(\omega\times w,v)+b(w_*,w,v)=<F,v>,\forall v\in\,V_m ,\eqno{(6.12)}$$
(6.12) is uniquely solvable because $w_*\in\,S_\rho$ because the
$w=0$ is only solution of the corresponding homogeneous equation
$(F=0)$. Indeed, if $w_*$ satisfies (6.10) and $w$ satisfies (6.12)
with the $F=0$, taking (6.4) and
$$2(\omega\times w,w)=0 \eqno{(6.13)}$$
 into account, and combing (3.19) and (6.11) we assert
$$\begin{array}{ll}
c(\Omega,\Theta)\mu\|\nabla w\|^2_{0,\Omega}\leq |b(w_*,w,w)|\leq
c_1(\Omega,K_0) \|w_*\|_{L^6}\|\nabla w\|_0\|w\|_2\leq c_1\|\nabla
w_*\|_0\|\nabla w\|^2_0\\
 \qquad\leq
c_1(\Omega,K_0)\frac{c(\Omega,\Theta)\mu}{2c_1}\|\nabla
w\|^2_0,\end{array}\eqno{(6.14)}$$ This implies that $w=0$. In order
to apply Brouwer's fixed point theorem, we have to show the mapping
$w_*\Rightarrow w$ take the ball $S_\rho$ defined by (6.10) into
itself, suppose that $w_*$ satisfies (6.10). Similarly to (6.13),
for nonhomogeneous equation $(F\neq 0)$ ,we obtain
$$\begin{array}{ll}
c(\Omega,\Theta)\mu\|\nabla w\|^2_{0,\Omega}\leq
|b(w_*,w,w)|+|<F,w>|\leq\\
\qquad c_1(\Omega,K_0)\|\nabla w_*\|_0\|\nabla
w\|^2_0+c_2(\Omega,K_0)\|F\|_*\|\nabla w\|_0,\end{array}$$
$$c(\Omega,\Theta)\mu\|\nabla w\|_{0,\Omega}\leq
 c_1(\Omega,K_0)\|\nabla w_*\|_0\|\nabla w\|_0+c_2(\Omega,K_0)\|F\|_*,\eqno{(6.15)}$$
 Here note that
$$Fw=g_{ij}F^iw^j=(\delta_{\alpha\beta}+r^2\Theta_\alpha\Theta_\beta)F^\alpha
w^\beta+\varepsilon\Theta_\alpha(F^\alpha
w^3+F^3w^\alpha)+\varepsilon^2r^2F^3w^3$$

 Therefore
 $$\begin{array}{ll}
 \|\nabla w\|_0\leq\frac{c_2\|F\|_*}{c(\Omega,\Theta)\mu-c_1\|\nabla w\|_0}\leq
 \frac{c_2\|F\|_*}{c(\Omega,\Theta)\mu-\frac{c(\Omega,\Theta)\mu}{2}
 [1-\sqrt{(1-4c_1c_2\|F\|_*/(c(\Omega,\Theta)\mu)^2)}]}\\
\quad
=\frac{c(\Omega,\Theta)\mu}{2c_1}[1-\sqrt{(1-\frac{4c_1c_2}{(c(\Omega,\Theta)\mu)^2}\|F\|_*)}]
 \end{array}$$
 Thus Brouwer's fixed point theorem can be applied and gives the
 existence of Galerkin approximations satisfying (6.8). Hence by a
 standard compactness argument there is at least a subsequence of
 the Galerkin approximation converging to a weak solution
 $w\in\,V(\Omega)$ of the steady problem (3.11):
 $$\left\{\begin{array}{ll}
\mbox{Find}\,(w,p),w\,\in\,V(\Omega),p\in\,L^2(\Omega),\mbox{such that}\\
a(w,v)+2(\omega\times w,v)+b(w,w,v)+\\
\quad\quad -(p, \div v)=<F,v>,\quad
\forall\,v\,\in\,V(\Omega),\\
(q,\div w)=0,\quad \forall\,q\,\in\,L^2(\Omega),\\
\end{array}\right.$$
   Its smoothness is easily proven if one obtain a further estimate
   from the Galerkin approximations by setting $v=-Aw$ in (6.9).
   This gives
   $$\mu\|Aw\|_0^2=-2(\omega\times w,Aw)-b(w,w,Aw)+<F,Aw>,\eqno{(6.14)}$$
   Because  $Aw$ is solenoidal, one has the rather unusual trace
   estimate
   $$|2(\omega\times w,Aw)|\leq c_3\|w\|_0\|Aw\|_0,\quad
   |<F,Aw>|\leq c_3\|F\|_*\|Aw\|_0,\eqno{(6.15)}$$
   which we combine with (6.14) and Agmon inequality
   $$\|w\|_{\infty}\leq c_4\|\nabla
   w\|_0^{1/2}\|Aw\|_0^{1/2},\quad\forall w\in\,D(A),\eqno{(6.16)}$$
   to get
   $$\mu\|Aw\|_0^2\leq c_4\|\nabla
   w\|_0^{3/2}\|Aw\|_0^{3/2}+c_3\|w\|_0\|Aw\|_0+c_3\|F\|_*\|Aw\|_0.\eqno{(6.17)}$$
   Then, by using Young's inequality, we obtain
   $$\mu\|Aw\|_0\leq \frac{2c^2_4}{\mu^2}\|\nabla
   w\|_0^3+\frac{8c^2_3}{\mu}\|w\|_0^2+\frac{8c^2_3}{\mu}\|F\|^2_*.\eqno{(6.18)}$$
which is then inherited by the solution. The full classical
smoothness of the solution can now be obtained using the
$L^2$-regularity theory for the steady Stokes equations. This
completes the proof of Theorem \ref{thm9}.
\end{proof}

It is clear that the bound in (6.10) is not uniform with respect to
$\Theta$. Next our gaol is to prove the the solution $w(\Theta)$ of
Navier-Stokes equations (3.11) with mixed boundary condition  is
uniformly bounded. At first we prove

\begin{lem}\label{lem8}
Assume that
$$\Theta\in\,{\cal
F}=\{\phi\in\,C^1(\Omega),\sup\limits_D|\nabla\phi|\leq
K_1:=\frac1{2\sqrt{5}r}, \}\eqno{(6.19)}$$ then
$$\begin{array}{ll}
a(w(\Theta),w(\Theta))\geq
\mu/2(2e_{\alpha\beta}(w)e_{\alpha\beta}(w)+e_{3\alpha}(w)e_{3\alpha}(w)+e_{33}(w)e_{33}(w))
,\quad\forall w\in V(\Omega)
\end{array}\eqno{(6.20)}$$
where $a(\cdot,\cdot)$ is defined by (6.1) and $\|\cdot\|_\Omega$ is
defined by (3.25).
\end{lem}
\begin{proof}
 By virtue of (4.6)(2.1)$(\lambda=0)$ we claim
$$\begin{array}{ll}
a(w,w)&=2\mu[e_{\alpha\beta}(w)e_{\alpha\beta}(w)+(2\varepsilon^{-2}\Theta_\lambda\Theta_\sigma+2\delta_{\lambda\sigma}
g^{33})e_{3\lambda}(w)e_{3\sigma}(w)+g^{33}g^{33}e_{33}(w)e_{33}(w)\\
&-4\varepsilon^{-1}\Theta_\alpha
e_{3\lambda}(w)e_{\alpha\lambda}(w)+2\varepsilon^{-2}\Theta_\lambda\Theta_\sigma
e_{\lambda\sigma}(w)e_{33}(w)-4\varepsilon^{-1}\Theta_\lambda
e_{3\lambda}(w)e_{33}(w)],
\end{array}\eqno{(6.21)}$$
By Cauchy inequality and Yang inequality we obtain
$$\begin{array}{ll}
4\varepsilon^{-1}\Theta_\lambda e_{3\lambda}(w)e_{33}(w)&\leq
8\varepsilon^{-2}(\Theta_\lambda
e_{3\lambda}(w))^2+\frac12(g^{33}e_{33}(w))^2,\\
4\varepsilon^{-1}\Theta_\alpha
e_{3\lambda}(w)e_{\alpha\lambda}(w)&\leq
4\sqrt{\varepsilon^{-2}\sum\limits_{\alpha\lambda}(\Theta_\alpha
e_{3\lambda}(w))^2}\sqrt{\sum\limits_{\alpha\lambda}(e_{\alpha\lambda}(w))^2}\\
&\leq
16\varepsilon^{-2}|\nabla\Theta|^2(e_{3\lambda}(w)e_{3\lambda}(w))+\frac14e_{\alpha\lambda}(w)e_{\alpha\lambda}(w),\\
2\varepsilon^{-2}\Theta_\lambda\Theta_\sigma
e_{\lambda\sigma}(w)e_{33}(w)&\leq
2\sqrt{\varepsilon^{-4}\sum\limits_{\lambda,\sigma}(\Theta_\lambda\Theta_\sigma
e_{33}(w))^2}\sqrt{\sum\limits_{\lambda,\sigma}(e_{\lambda\sigma}(w))^2}\\
&\leq
4\varepsilon^{-4}|\nabla\Theta|^4e_{33}(w)e_{33}(w)+\frac14e_{\lambda\alpha}(w)e_{\alpha\lambda}(w),
\end{array}\eqno{(6.22)}$$
To sum up and thanks
$g^{33}=\varepsilon^{-2}r^{-2}+\varepsilon^{-2}|\nabla\Theta|^2$ we
assert that
$$\begin{array}{ll}
a(w,w)&\geq
2\mu[\frac12e_{\alpha\beta}(w)e_{\alpha\beta}(w)+[-6\varepsilon^{-2}\Theta_\lambda\Theta_\sigma+
2\varepsilon^{-2}r^{-2}\delta_{\lambda\sigma}-14\varepsilon^{-2}|\nabla\Theta|^2\delta_{\lambda\sigma}]
e_{3\lambda}(w)e_{3\sigma}(w)\\
&+(\frac12g^{33}g^{33}-4\varepsilon^{-4}|\nabla\Theta|^4)e_{33}(w)e_{33}(w)],\\
\end{array}$$ Note
$$\Theta_\lambda e_{3\lambda}(w)\Theta_\sigma e_{3\sigma}(w)\leq
(\sqrt{\sum\limits_{\lambda}\Theta_\lambda^2}\sqrt{\sum\limits_{\lambda}e^2_{3\lambda}(w)})^2
\leq |\nabla\Theta|^2e_{3\lambda}(w)e_{3\lambda}(w)$$
$$(\frac12g^{33}g^{33}-4\varepsilon^{-4}|\nabla\Theta|^4)=\frac{\varepsilon^{-4}}{2}(
r^{-4}+2r^{-2}|\nabla\Theta|^2-7|\nabla|^4)$$ Therefore
$$\begin{array}{ll} a(w,w)&\geq
2\mu[\frac12e_{\alpha\beta}(w)e_{\alpha\beta}(w)+
2\varepsilon^{-2}r^{-2}[1-10r^2|\nabla\Theta|^2]
e_{3\lambda}(w)e_{3\lambda}(w)\\
&+\frac{\varepsilon^{-4}r^{-4}}{2}(
1+2r^{2}|\nabla\Theta|^2-7r^4|\nabla\Theta|^4)e_{33}(w)e_{33}(w)]\\
\end{array}$$
It is obvious that
$$\begin{array}{ll}
2\varepsilon^{-2}r^{-2}[1-10r^2|\nabla\Theta|^2]\geq
\frac14\varepsilon^{-2}r^{-2},\quad\mbox{if}\, |\nabla\Theta|\leq
\frac1{2\sqrt{5}r};\\
 \frac{\varepsilon^{-4}r^{-4}}{2}(
1+2r^{2}|\nabla\Theta|^2-7r^4|\nabla\Theta|^4)\geq
\frac14\varepsilon^{-4}r^{-4},\quad\mbox{if}\,|\nabla\Theta|\leq\frac2{7r},
\end{array}$$ Hence If $\Theta\in {\cal F} $ then
$$\begin{array}{ll} a(w,w)&\geq

2\mu[\frac12e_{\alpha\beta}(w)e_{\alpha\beta}(w)+\frac14e_{3\lambda}(w)e_{3\lambda}(w)
+\frac14e_{33}(w)e_{33}(w)],
\end{array}$$ End the proof.
\end{proof}

\begin{lem}\label{lem9}
Assumption in Lemma \ref{lem4} is held such that the constant
$K_0$ satisties
$$K_0\leq \frac1{\sqrt{2\beta_0(\Omega)}},\eqno{(6.23)}$$
where $\beta_0(\Omega)$ is a constants depending $\Omega$ only .Then
we have following unform coerciveness for bilinear form
$a(\cdot,\cdot)$
$$a(w,w)\geq \frac{\mu}{4}\|w\|^2_\Omega\geq \frac{\mu}{4}C_1(\Omega)\|w\|^2_{1,\Omega},
\quad\forall\,w\in\,H^1(\Omega)^3,\eqno{(6.24)}$$ where
$\alpha_0(\Omega)$ is a constant independent of $\Theta$ and
$\|w\|^2_\Omega$ is defined by (3.25).
\end{lem}

\begin{proof}
In view of (3.16),
$$\begin{array}{ll}
\|e_{\alpha\beta}(w)\|^2_{0,\Omega}&=\|\varphi_{\alpha\beta}(w)\|^2_{0,\Omega}
+\|\psi^\lambda_{\alpha\beta}(w)\Theta_\lambda\|^2_{0,\Omega}+
\|\psi^{\lambda\sigma}_{\alpha\beta}(w)\Theta_\lambda\Theta_\sigma\|^2_{0,\Omega}
+\|e^*_{\alpha\beta}(w)\|^2_{0,\Omega}\\
&+2(\varphi_{\alpha\beta}(w),\psi^\lambda_{\alpha\beta}(w)\Theta_\lambda)
+2(\varphi_{\alpha\beta}(w),\psi^{\lambda\sigma}_{\alpha\beta}(w)\Theta_\lambda\Theta_\sigma)\\
&+2(\varphi^{\nu\mu}_{\alpha\beta}(w)\Theta_\nu\Theta\mu,\psi^\lambda_{\alpha\beta}(w)\Theta_\lambda)
+2(\widetilde{e}_{\alpha\beta}(w),e^*_{\alpha\beta}(w)),
\end{array}\eqno{(6.25)}$$
Using (3.18) it is easy to verify
$$\begin{array}{ll}
\|\psi^\lambda_{\alpha\beta}(w)\|^2_{0,\Omega}\leq
\beta_0(\Omega)\|\varphi_{\alpha3}(w)\|^2_{0,\Omega},\quad
\|\psi^{\lambda\sigma}_{\alpha\beta}(w)\|^2_{0,\Omega}\leq
\beta_0(\Omega)(\|\varphi_{\alpha\beta}(w)\|^2_{0,\Omega}+\|\varphi_{33}(w)\|^2_{0,\Omega}),\\
\|\psi^\lambda_{3\alpha}(w)\|^2_{0,\Omega}\leq
\beta_0(\Omega)(\|\varphi_{\alpha\beta}(w)\|^2_{0,\Omega}+\|\varphi_{33}(w)\|^2_{0,\Omega})\\
\|\psi^{\lambda\sigma}_{3\alpha}(w)\|^2_{0,\Omega}\leq\beta_0(\Omega)\|\varphi_{\alpha3}(w)\|^2_{0,\Omega},\quad
\|\psi^\lambda_{33}(w)\|^2_{0,\Omega}\leq\beta_0\|\varphi_{3\alpha}(w)\|^2_{0,\Omega},\\
\|e^*_{\alpha\beta}(w)\|^2_{0,\Omega}\leq\beta_0(\Omega)K_0^2\|\varphi_{\alpha\beta}(w)\|_{0,\Omega},\quad
|(\widetilde{e}_{\alpha\beta}(w), e^*_{\alpha\beta}(w))|\leq
K_0^{1/2}\beta_0(\Omega)\|w\|^2_D,
\end{array}\eqno{(6.26)}$$

\begin{remark} the constant $\beta_0(\Omega)$ represents different meaning at
different place. \end{remark}

 By means of Schwarcs inequality, from
(6.25)(6.26) we claim
$$\|e_{\alpha\beta}(w)\|^2_{0,\Omega}\geq
\|\varphi_{\alpha\beta}(w)\|^2_{0,\Omega}-
K^{1/2}_0\beta_0(\Omega)\|w\|_\Omega^2,\eqno{(6.27)}$$ By similar
manner we imply
$$\|e_{3\alpha}(w)\|^2\geq\|\varphi_{3\alpha}(w)\|^2_{0,\Omega}-K_0^{1/2}\beta_0(\Omega)\|w\|^2_\Omega,\quad
\|e_{33}(w)\|^2\geq\|\varphi_{33}(w)\|^2_{0,\Omega}-K_0\beta_0(\Omega)\|w\|^2_\Omega,$$
Let return to (6.20), it yields that
$$\begin{array}{ll}
a(w,w)\geq
\frac{\mu}{2}[\|w\|^2_{0,\Omega}-K^{1/2}_0\beta_0(\Omega)\|w\|^2_{0,\Omega}].$$
\end{array}$$
This yields that if (6.23) is satisfied then ,by Lemma \ref{lem3},
it yields (6.24) .
\end{proof}

From Theorem \ref{thm6} and Lemma \ref{lem5} we conclude

\begin{thm}\label{thm10}Assume the the assumptions in Theorem \ref{lem6} are satisfied and
(6.23) held. The $w(\Theta)$ is a solution of rotating Navier-Stokes
equations (3.11) associated with $\Theta\in C^2(\Omega,K_0)$. Then
the following estimate  held
$$\| w(\Theta)\|_{\Omega}\leq
\frac{\mu}{2c_1}[1-\sqrt{(1-\frac{4c_1c_2\|F\|_*}{(\mu)^2})}],\eqno{(6.28)}$$
The solution is sequently weak continuous with respect to $\Theta$,
i.e. if sequence $ \{\Theta_i\}$ is weak convergent in $H^1(D)$ then
there exist convergent subsequent in $\{w(\Theta_i)\}$ in
$H^1(\Omega)$.
\end{thm}
%%%%%%%%%%%%%%%%%%%%%%%%%%

\end{document}